\documentclass[review]{elsarticle}

\usepackage[english]{babel}
\usepackage[show]{chato-notes}
\usepackage{mathtools}
\usepackage{amsmath}
\usepackage{amssymb}
\usepackage{hyperref}
\usepackage{subfigure}
\usepackage{tikz}
\usetikzlibrary{backgrounds}
\usetikzlibrary{arrows.meta}
\usetikzlibrary{shapes,snakes}

\newtheorem{theorem}{Theorem}

\newtheorem{remark}[theorem]{Remark}

\begin{document}

\begin{frontmatter}

\title{Structural Systems Theory: an overview of the last 15 years\tnoteref{mytitlenote}}
\tnotetext[mytitlenote]{}

\author{Guilherme Ramos\fnref{myfootnote}}
\address{Department of Electrical and Computer Engineering, Faculty of Engineering, University of Porto, Portugal}
\author{A. Pedro Aguiar\fnref{myfootnote2}}
\address{Department of Electrical and Computer Engineering, Faculty of Engineering, University of Porto, Portugal}
\author{S\'ergio Pequito}
\address{Delft Center for Systems and Control, Delft  University of Technology, 2600 AA Delft, The Netherlands}

\begin{abstract}
In this paper, we provide an overview of the research conducted in the context of structural systems since the latest survey by Dion \emph{et al.} in 2003. We systematically consider all the papers that cite this survey as well as the seminal work in this field that took place on and after the publication of the later survey, are published in peer-reviewed venues and in English.

Structural systems theory deals with parametric systems where parameters might be unknown and, therefore, addresses the study of systems properties that depend only on the system’s structure (or topology) described by the inter-dependencies between state variables. Remarkably, structural systems properties hold generically (i.e., almost always) under the assumption that parameters are independent. Therefore, it constitutes an approach to assess necessary conditions that systems should satisfy. 

In recent years, structural systems theory was applied to design systems that attain such properties, as well as to endure resilient/security and privacy properties. Furthermore, structural systems theory enables the formulation of such topics as combinatorial optimization problems, which allow us to understand their computational complexity and find algorithms that are efficiently deployed in the context of large-scale systems. 

In particular, we present an overview of how structural systems theory has been used in the context of linear time-invariant systems, as well as other dynamical models, for which a brief description of the different problem statements and solutions approaches are presented. Next, we describe recent variants of structural systems theory, as well as different applications of the classical and new approaches. Finally, we provide an overview of recent and future directions in this field.
\end{abstract}


\end{frontmatter}


{\scriptsize
\tableofcontents
}

\section{Introduction}



Structural systems theory is nowadays equipped with powerful tools to perform analysis, design, and optimization of dynamical systems when their parameters cannot be considered accurate or even completely unknown, yet their inter-dependencies can be assumed to be known (i.e., the systems structure or topology), and under the assumption that those parameters are independent of each other.
Structural systems theory and the notion of structural systems properties were first coined by Lin in 1974~\cite{lin1974structural}, where the notion of structural controllability was introduced for single-input single-output systems. It was only in 1976, that Shields and Pearson extended the notion of structural controllability for multi-input multi-output systems and presented rigorous characterizations for it (i.e., necessary and sufficient conditions)~\cite{shields1976structural}.

A particularly interesting aspect of structural systems theory is the fact that systems' properties hold in general (i.e., for almost all possible realizations of the independent parameters).  It worth mention that controllability, and subsequently observability, are known to be properties that hold generically since at least 1962~\cite{markus1962existence}.  
Notwithstanding, this main advantage is also a limitation in that it only allows a qualitative assessment of systems properties. In other words, since there is no need to have the parameters' values, we cannot quantify some system properties (e.g., controllability/observability energy through the respective Grammians). Nonetheless, we argue that in many cases this limitation is relative, as such quantitative assessment  might be prohibitive when dealing with large scale systems -- see Remark~\ref{rem:1} and~\ref{rem:2}.

In recent years, structural systems theory was considered to design systems that attain such properties, as well as to endure resilient/security and privacy properties. Furthermore, structural systems theory enables the formulation of such problems as combinatorial optimization problems, which enable us to understand their computational complexity and find algorithms that are efficiently deployed in the context of large-scale systems. In this paper, we provide an overview on how structural systems theory has been used in the context of linear time-invariant systems, as well as other dynamical models, for which a brief description of the different problem statements and solutions approaches are presented. Next, we describe recent variants of structural systems theory, as well as different applications of the classical and new approaches. Finally, we provide an overview of recent and future directions in this field for those looking forward to initiate their research in the topic, or simply, to leverage tools provided by structural systems towards deriving \emph{non-}structural properties.\\

\noindent \textbf{Survey methodology:} In this overview, we consider all the papers written in English that make use or leverage structural systems theory, as well as their direct variants and applications.
In particular, we take into consideration all peer-reviewed papers in press that cited either Lin's paper~\cite{lin1974structural}, or Dion \emph{et al.}~\cite{dion2003generic}, that were available from 2003 to the end of March 2020. Also, we consider only the conference publications when these did not seem to have a published journal version. Besides, we do not cite dissertations as they contain parts of papers that are cited. Lastly, we would like to disclaim that the aim of this study is not to provide an in-depth overview of all results, but rather a (possibly bias) selection of results that could help further developments in the context of structural systems theory and its use to show \emph{non-}structural system properties. In the same line, we try to report the results found during the survey without judging their validity and/or usefulness.

\section{Linear Time-invariant Systems}\label{sec:LTI}

A variety of dynamical systems including mechanical and electrical systems, multi-agents, social networks, biological systems, and many others, can be modeled/described by a Linear Time-Invariant (LTI) system, whose state-space representation for the case of discrete-time is given by
\begin{equation}\label{eq:lti}
    \begin{split} x(t+1)&=Ax(t)+Bu(t),\\
        y(t)&=Cx(t)+D u(t), \ t=0,1,\ldots,
    \end{split}
\end{equation}
where $x(t)\in\mathbb R^n$ is the state, $x(0)=x_0$ is the initial condition, $u(t)\in\mathbb R^p$ is  the input, and $y(t)\in\mathbb R^q$ is the output. 
In what follows, we focus on the discrete-time LTI as most results readily apply to continuous-time~LTI. 
As such, we omit the reference to the fact that we are dealing with discrete-time, and we will emphasize when the results only apply to either discrete- or continuous-time~LTI.

In the context of structural systems theory, the goal is to consider only the \emph{structural pattern} of the model parameters in~\eqref{eq:lti} (i.e., the matrices $A$, $B$, and~$C$). 
Specifically, we only need to check  if the parameters are equal to zero -- referred to as \emph{fixed zeros}) -- and denoted by $0$, or if the parameter can take any value (including zero) -- referred to as a \emph{free} parameter and, with some abuse of terminology, a \emph{nonzero} parameter --  and denoted by $\star$. 
Recall that it is implicitly assumed that such nonzero parameters are independent of each other.

In what follows, when considering only the structural pattern, we will denote vectors and matrices with a bar on top (e.g., $\bar v$ and $\bar A$).

\begin{figure}[!h]
\centering
\includegraphics[width=0.6\textwidth]{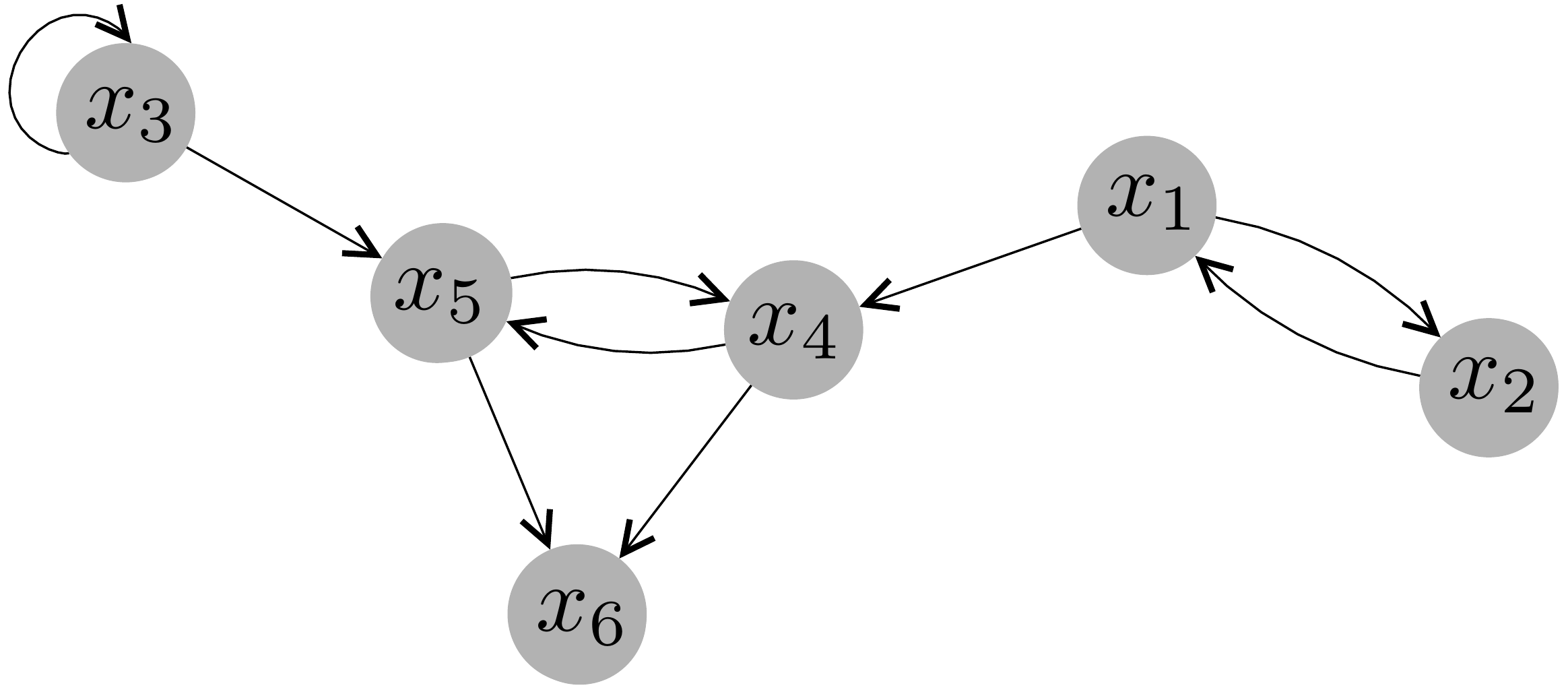}
\caption{Digraph representation of the structural matrix $\bar A$, $\mathcal G(\bar A)$.}
\label{fig:state}
\end{figure}

An appealing feature of dealing with structural systems and their structural patterns is through a natural (and equivalent) representation as a directed graph (digraph) that consists in a set of vertices (or nodes) and edges~\cite{cormen2009introduction}. 
The digraph representation of the systems dynamics is provided by the \emph{state digraph}  $\mathcal G(\bar A)=(\mathcal X,\mathcal E)$, where $\mathcal X=\{x_1,\ldots,x_n\}$ denote vertices labeled by the states  of the system, and $\mathcal E=\{(x_i,x_j)\,:\,\bar A_{ij}\neq 0\}$ are the edges capturing the  dependencies between the states.  
For example, suppose we have an LTI system with a dynamics matrix and input matrix structural pattern given by
\[
\bar A=\left[
\begin{array}{cccccc}
    0      & \star & 0     & 0     & 0     & 0 \\
    \star  & 0     & 0     & 0     & 0     & 0 \\
    0      & 0     & \star     & 0     & 0     & 0 \\
    \star  & 0     & 0     & 0     & \star & 0 \\
    0      & 0     & \star & \star & 0     & 0 \\
    0      & 0     & 0     & \star & \star & 0 \\
\end{array}
\right], \,
\bar B =
\left[
\begin{array}{cc}
    0 & \star \\
    0 & 0 \\
    \star & 0 \\
    0 & 0 \\
    0 & 0 \\
    0 & 0
\end{array}
\right]
\] 
and, in Fig.~\ref{fig:state}, we depict the state digraph  $\mathcal G(\bar A)$. 
In the input matrix $\bar B$, there are two inputs, the first one, $u_1$, is assigned to state variable $x_3$ and the second one, $u_2$, is assigned to state variable $x_1$. 
Analogously, we may consider the structural pattern of the output matrix, where again an entry is non-zero ($\star$) if in the output matrix ($C$) the entry is nonzero and an entry is zero otherwise. 

\emph{Input and output digraphs} can be similarly obtained by considering the input and output matrices, respectively. 
In these cases, the dependencies between the input  and state vertices, as well as state and output vertices are unidirectional, recall the matrix $\bar B$ in the previous example. 
Specifically, the input vertices only have outgoing edges, whereas the output vertices only have incoming edges. In other words, the edges from the inputs to the states represent which states are under direct influence of the actuators/controller, whereas the edges to the output vertices encode the observations/measurements of the states obtained by the sensors. 
Therefore, the digraph representation yields a simple way to visualize the dependencies between the variables of the system. Furthermore, it enable us to study (structural) systems properties as graphical properties for which there might exist efficient algorithms to assess them as well as to design systems that have such properties. 

\subsection{Structural controllability} The first structural property that one can encounter in structural systems theory is that of \emph{structural controllability}. 
Structural controllability of the system described by the dynamics and input matrices with structural pattern $(\bar A,\bar B)$ is attained if and only if there exists a controllable system described by $(A,B)$ with the same structural pattern as  $(\bar A,\bar B)$. 
As previously mentioned, we can invoke measure theoretical arguments to establish that structural controllability is a generic property. 
In other words, if there is a pair $(A,B)$ with the same structural pattern as $(\bar A,\bar B)$, then almost all pairs $(A',B')$ with the same structural pattern lead to controllable systems~\cite{lin1974structural,dion2003generic}. 
In particular, it readily follows that if one replaces the nonzero entries in $(\bar A,\bar B)$ at random by entries drawn from continuous distributions, such pair will be almost surely controllable. 
Note also that if a pair $(\bar A,\bar B)$ is not structurally controllable, then there exists no realization of such a structural pattern that yields a controllable pair. 
Thus, structural controllability entails a necessary condition for controllability.

\begin{figure}[!h]
\centering
\includegraphics[width=0.35\textwidth]{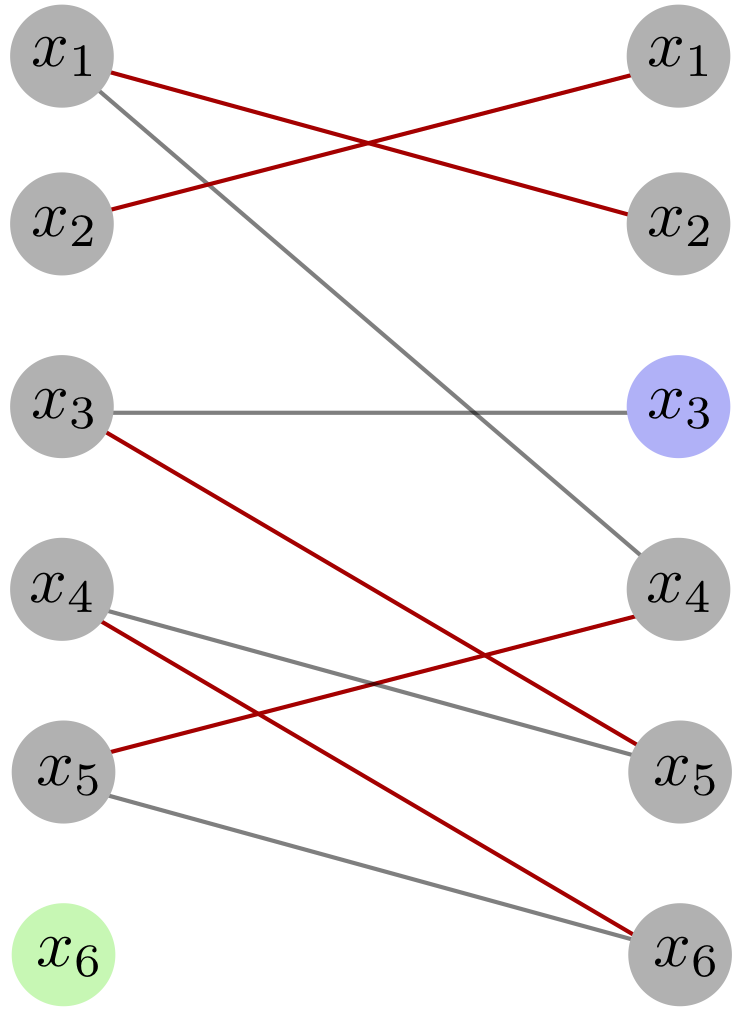}
\caption{Bipartite graph representation of $\bar A$, $\mathcal B(\bar A)$. The edges in red corresponds to one (from possible multiple choices) of a maximum matching of $\mathcal B(\bar A)$, the vertex in blue (on the right hand-side) is a right-unmatched vertex and the vertex in green (on the left hand-side) is a left-unmatched vertex. }
\label{fig:bip}
\end{figure}

To assess structural controllability, one can use a variety of graph theoretical tools like path and cycle decompositions and the use of an auxiliary \emph{(undirected) bipartite graph} on which we can find \emph{maximum matchings}~\cite{dion2003generic}. 
The bipartite graph representation associated with the state digraph $\mathcal G(\bar A)$ is denoted by $\mathcal B(\bar A)=(\mathcal X^l,\mathcal X^r,\mathcal E)$, where two copies of the state vertices ($\mathcal X^l,\mathcal X^r$) are considered to represent those state vertices in one set that can be visually presented on the left (denoted by $\mathcal X^l$), and another on the right (denoted by $\mathcal X^r$) and the set of undirected edges as the ordered pairs $\mathcal E=\{(x_i,x_j)\,:\,\bar A_{ij}\neq 0\}$; see Fig.~\ref{fig:bip}. With some abuse of terminology, to avoid over-notation, the bipartite graph is often presented as $\mathcal B(\bar A)=(\mathcal X,\mathcal X,\mathcal E)$, where left and right sets of vertices are implicitly indicated by order of appearance in the tuple that defines the bipartite graph.

A maximum matching on a bipartite graph corresponds to the most extensive collection of edges that do not share vertices within each vertex sets -- see the edges depicted in red in Fig.~\ref{fig:bip}. 
The set of left vertices that are not present at the edges that belong to the maximum matching are called left-unmatched vertices. 
Analogously, the set of right-unmatched vertices is the set of right vertices are not present in the edges of the maximum matching. 
For instance, in Fig.~\ref{fig:bip}, the left-unmatched vertex $x_6$ is colored with green and the right-unmatched vertex $x_3$ is colored in blue. These concepts can be easily extended to the concepts of input and output bipartite graphs, with appropriate changes on the vertices set due to the existence of input and output
vertices on the input and output digraphs, respectively.

As previously mentioned, the role of the bipartite graph, and in particular of the maximum matching, is that of serving as a proxy for properties that can be stated on the state/input/output digraphs. In particular, a maximum matching enables us to identify a  decomposition into disjoint paths and cycles in the digraph, where the directed edges belonging to these are determined by the undirected edges which origin is on the left vertices set and the end on the right vertices set.
In fact, such decomposition has the minimum number of paths (possibly degenerated, i.e., single vertices) and an arbitrary number of cycles. 
Such decomposition in paths and cycles of the digraph plays a critical role in the notion of structural controllability~\cite{dion2003generic}. More generally, all maximum matchings can be described through the so-called Dulmage–Mendelsohn decomposition that plays a role in characterizing the fixed controllable subspace~\cite{commault2017fixed}. 
Nonetheless, it is possible to assess structural controllability by resorting to other methods, e.g., dynamic graph properties~\cite{van2019dynamic}.

At this point, the essential message to retain is that the structural system's properties can be easily stated in terms of graph-theoretical properties. However, the representation choice leads to different (discrete) combinatorial optimization problems -- some of which can be efficiently solved whereas others enable us only to obtain efficient approximate solutions with possibly optimality guarantees.  Whereas the computational efficiency of different representation choices for assessing structural systems properties may not differ much, they can make the difference when designing systems that possess such structural properties. To make a parallel that is more familiar for this paper's possible audience, a (classical) optimization problem may seem nonconvex, yet it might admit a convex representation that ensures the use of  computationally efficient numerical algorithms to determine its solution. Nonetheless, different convex representations might be possible for the same problem (i.e., linear program, quadratic program, geometric, second-order cone program, or semi-definite program), yet they require numerical algorithms that have different computational efficiency. 

Additionally, as systems increase their dimension, it is of interest to determine  conditions that allow to assess structural controllability efficiently by possibly considering composite systems or leveraging distributed algorithms~\cite{carvalho2017composability}.

In what follows, we will see that computational complexity theory enables us to quantify how computationally efficient we can find a solution to different problems.

\subsubsection{Computational complexity}\label{sub:complexity}

In computational complexity theory, we define a \emph{complexity class} as a set of computational problems of similar resource-based complexity~\cite{arora2009computational}. 
The two most commonly analyzed resources are time and memory. 

A usual notation in computational complexity theory is the big-$\mathcal O$ notation. 
Let $f$ and $g$ be two function with domain $\mathbb Z^+$. We say that $f(n)\in\mathcal O(g(n))$ whenever there is a constant number $k\in\mathbb Z^+$ and an order $n_0\in\mathbb Z^+$ such that for any $n\geq n_0$ it follows that $|f(n)|\leq |k g(n)|$. 
Hence, we may immediately identify some complexity classes (that can be in terms of time or memory):  
the \emph{constant class} is $\mathcal O(1)$; \emph{logarithmic class} (\textsc{Log}) -- $\mathcal O(\log n)$; \emph{linear class} (\textsc{Linear}) -- $\mathcal O(n)$; \emph{quadratic class} -- $\mathcal O(n^2)$; \emph{polynomial class} (\textsc{P}) -- $\mathcal O(n^c)$ for some $c\in\mathbb Z^+ $; and \emph{exponential class} (\textsc{Exp}) -- $\mathcal O(c^n)$ for some $c>1$. 
Another important time-complexity class is the \textsc{NP} class. A computational problem is in \textsc{NP} if there is no known polynomial algorithm to solve it. 

Observe that we cannot use more memory than the time we use. In other words, a computational problem that belongs to a particular time-complexity class, $\mathcal O(g(n))$, also belongs, in the worst case, to the memory-complexity class $\mathcal O(g(n))$, since in the number of time steps needed to solve the problem we can, at most, use the same number of memory. 
Further, we say that a computational problem can be solved \emph{efficiently} if it belongs to the polynomial class. 

\begin{figure}[!h]
    \centering
    \includegraphics[width=0.5\textwidth]{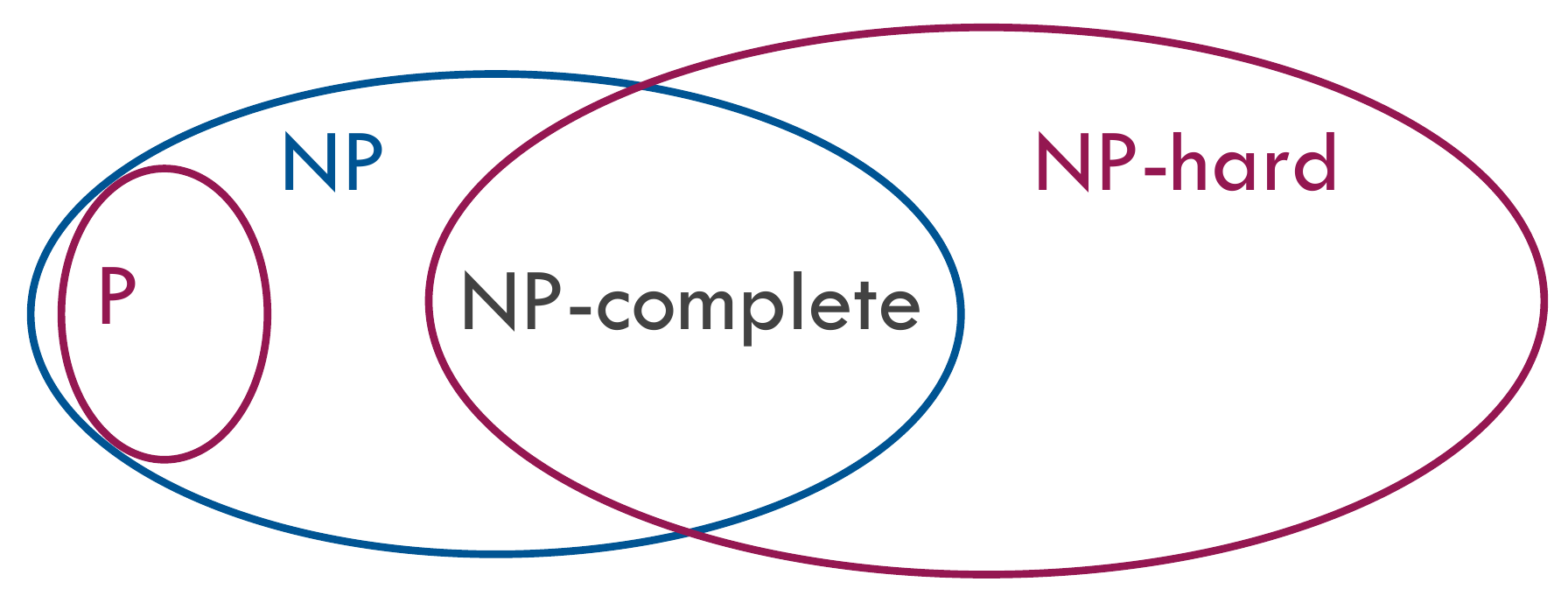}
    \caption{Schematic representation, assuming P$\neq$NP, of the relations between the complexity classes: P, NP, NP-hard and NP-complete.}
    \label{fig:complex_jungle}
\end{figure}

A computational problem $\mathcal P$ is reducible in polynomial-time to another $\mathcal P'$, $\mathcal P\prec_P \mathcal P'$, if there exists a procedure to transform the former to the latter using a polynomial number of operations on the size of its inputs. Such types of reductions are useful to determine the complexity class for which a particular problem belongs to~\cite{arora2009computational}. 
A problem $\mathcal P$ is in NP if, given a candidate solution to the problem, it can be verified if it is indeed a solution in polynomial time. In other words, the problems in NP are those whose solutions can be verified efficiently. 
A problem $\mathcal P$ is NP-hard when every problem $\mathcal P'$ in NP can be reduced in polynomial time to $\mathcal P$. 
Finally, a problem is NP-complete if it is both in NP and NP-hard. 
Alternatively, a problem $\mathcal P$ is NP-hard if there is an NP-complete problem $\mathcal P'$ such that $\mathcal P'\prec_P\mathcal P$. Assuming P$\neq$NP, in Fig.~\ref{fig:complex_jungle}, we depict the previous complexity classes and their relations. 

\subsection{Design systems to attain structural controllability}\label{sec:systemDesignsC}

In the context of designing structural systems that attain structural controllability, there are two main classes of problems: \emph{(i)} actuator placement (or
input/actuator selection) problem; and \emph{(ii)} dynamics topology design problem.



To elucidate the goals of the objectives and constraints in structural systems concerning \emph{non-}structural systems properties (e.g.,  controllability) consider the following problem statement: 
let $\Theta\equiv\Theta(A,B)$ be a collection of the interesting parameters in $A$ and/or $B$ and $f_\Theta$ be an objective function 
that assigns a cost to each combination of interesting parameters in $\Theta$, $f_\Theta:\mathbb{R}^{n\times n}\times \mathbb{R}^{n\times p}\rightarrow \mathbb{R}$ parameterized by $\Theta$, and where $(A,B)$ describes the model~\eqref{eq:lti}, 
we seek to solve problems of the form 
\begin{equation}\label{eq:prob_min_control}
\begin{split}
    \mathop{\min}_{\Theta}  & \quad f_\Theta(A,B) \\
\text{s.t. }& rank(\mathcal C(A,B)) =n,
\end{split}
\end{equation}
where the \emph{partial controllability matrix} is given by
\begin{equation}
\mathcal C_i(A,B)=[B\; AB\; A^2B\;\ldots\;A^{i-1}B], \quad i\in\mathbb{N}, 
\label{eq:partialcontrolmatrix}
\end{equation}
 and, subsequently, $\mathcal C(A,B)=\mathcal C_n(A,B)$ denotes the \emph{controllability matrix} used to establish the controllability of a system described by $(A,B)$ (i.e., the rank of the controllaiblity matrix equals the dimension $n$ of the state).

A questions that emerges is: ``Are there advantages of solving a version of  problem~\eqref{eq:prob_min_control} that considers only the structure, instead of solving problem~\eqref{eq:prob_min_control}?''

To answer the question, we start by noticing that this formulation has the following drawbacks. 

\begin{figure}[!h]
\subfigure[]{\label{fig:exp_motiv_1a}\includegraphics[width=0.47\textwidth]{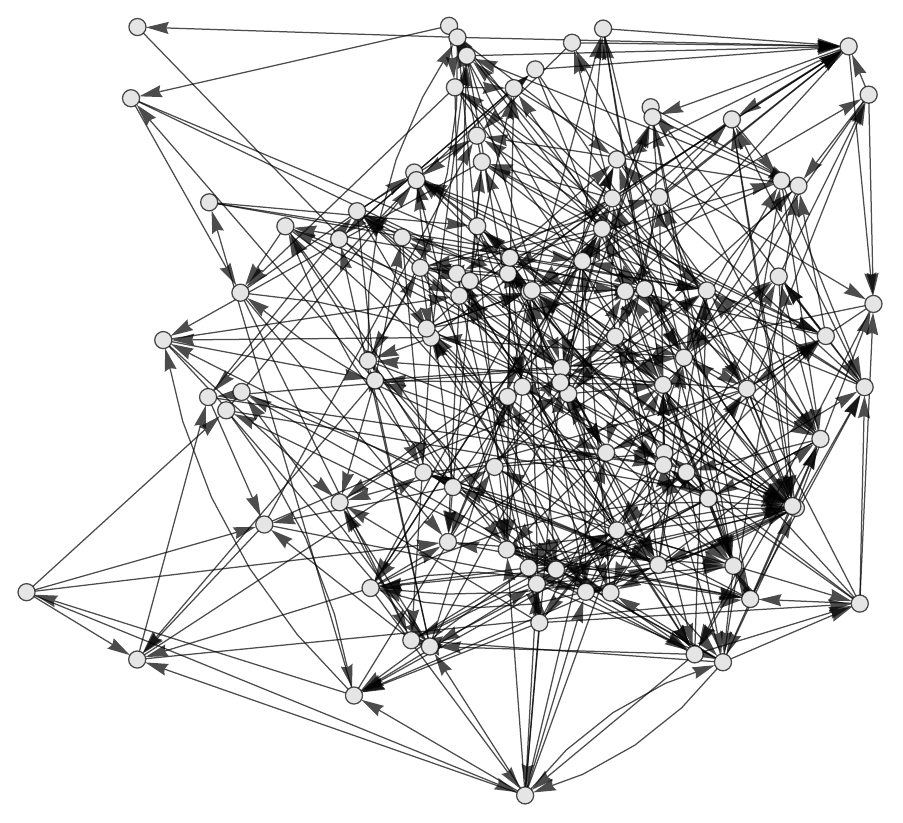}}
\hfill
\subfigure[]{\label{fig:exp_motiv_1b}\includegraphics[width=0.47\textwidth]{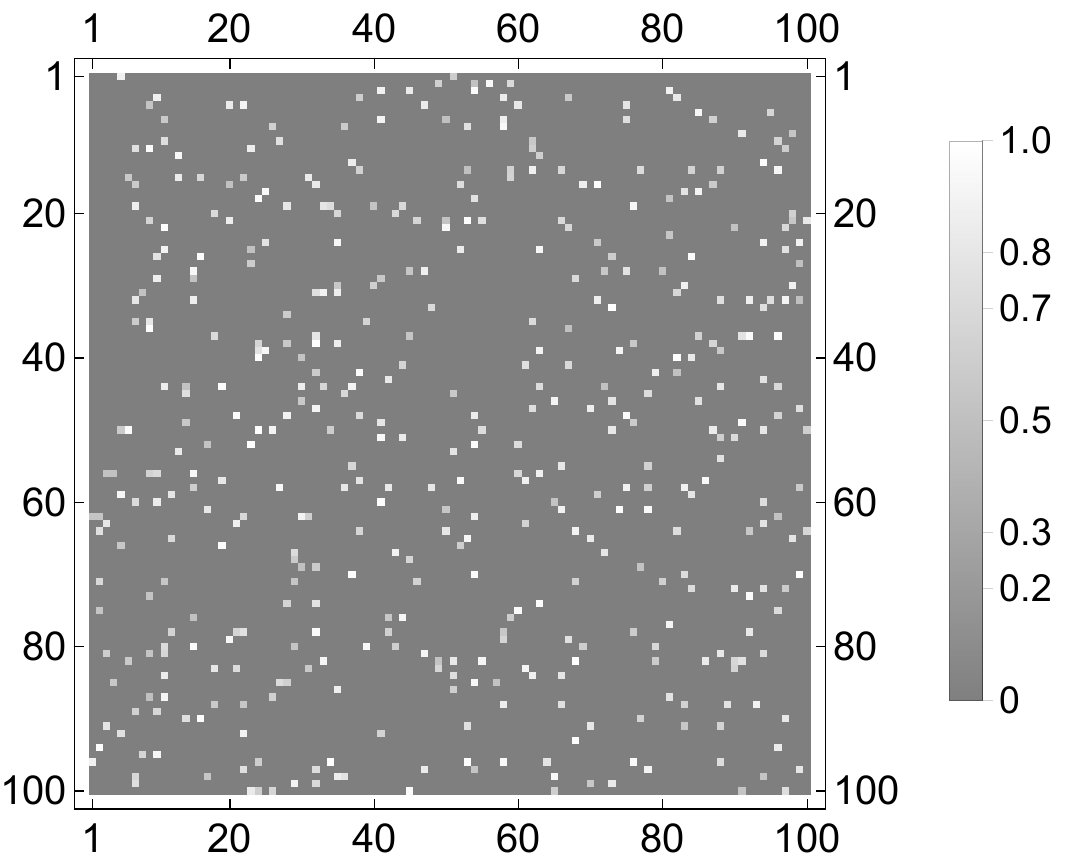}}
\caption{Randomly generated digraph with 100 vertices, in~(a), and repective adjacency matrix plot, in~(b). We can notice that the adjacency matrix is sparse. }
\label{fig:exp_motiv_1}
\end{figure}

\vspace{2mm}
\setlength{\fboxsep}{0pt}
\fcolorbox{white}[HTML]{f3f3f3}{\parbox{.94\textwidth}{
\begin{remark}[Controllability matrix and floating-point errors]\label{rem:1}
To illustrate some caveats that emerge when studying LTI systems' controllability properties. 
We start by exploring the controllability matrix's use as a controllability criterion for LTI systems. 
We generated the random and strongly connected digraph with 100 vertices, $\mathcal G_{100}\equiv \mathcal G(\bar A_{100})$, depicted in Fig.~\ref{fig:exp_motiv_1a}, whose adjacency matrix plot is given in Fig.~\ref{fig:exp_motiv_1b}. Now, consider  a matrix $A$ with the same sparseness as $\bar A$ such that $A_{ij}=1$ if $\bar A_{ij}\neq \star$ and $A_{ij}=0$ otherwise. 
Additionally, consider as the input matrix, the $100\times 1$ matrix $B$ with $B_{77}=1$ and $B_i=0$ for $i\neq 77$. The pair $(\bar A,\bar B)$ is structurally controllable. 
Under the described setup, if $A$ is the adjacency matrix of the digraph (i.e., the nonzero entries equal to 1), then we have that $rank(\mathcal C(A,B))=100$, yielding a controllable system. 
Subsequently, we add uniform noise to the nonzero parameters of $A$, uniform in the interval $]-10^{-6},10^{-6}[$, obtaining a matrix $A'$ also with the same sparsity of $\bar A$. 
Utilizing the same input matrix $B$, we compute the controllability matrix rank and obtain that $rank(\mathcal C(A',B))=1$. 
However, the rank should also be 100. 
This enormous gap between the controllability matrices is due to floating-point error propagation when computing the powers of matrix $A'$. 
Consequently, for large-scale systems, the controllability matrix $rank$ condition becomes unfeasible to be tested. 
Lastly, notice that the matrices $A$ and $A'$ are sparse and that more dense matrices would results in further floating-point error propagation. 
\end{remark}
}}

\begin{figure}[!h]
\subfigure[]{\label{fig:a}\includegraphics[width=0.35\textwidth]{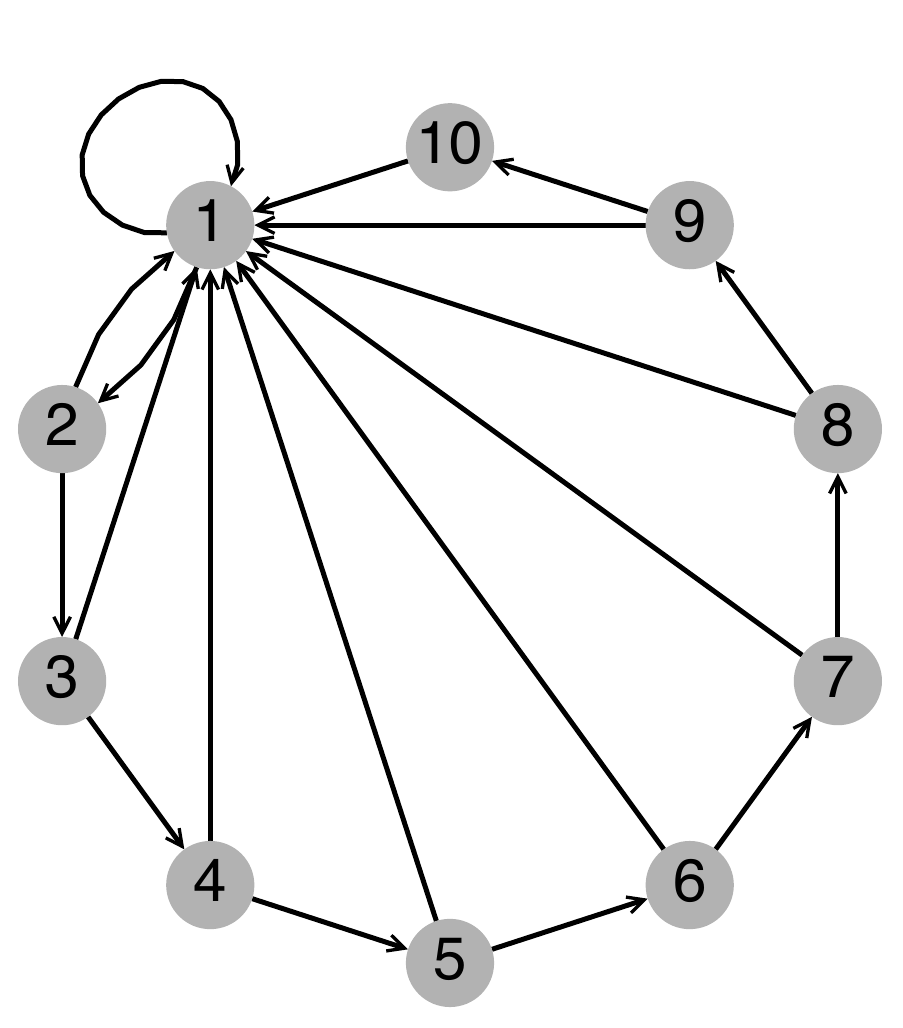}}
\hfill
\subfigure[]{\label{fig:a}\includegraphics[width=0.58\textwidth]{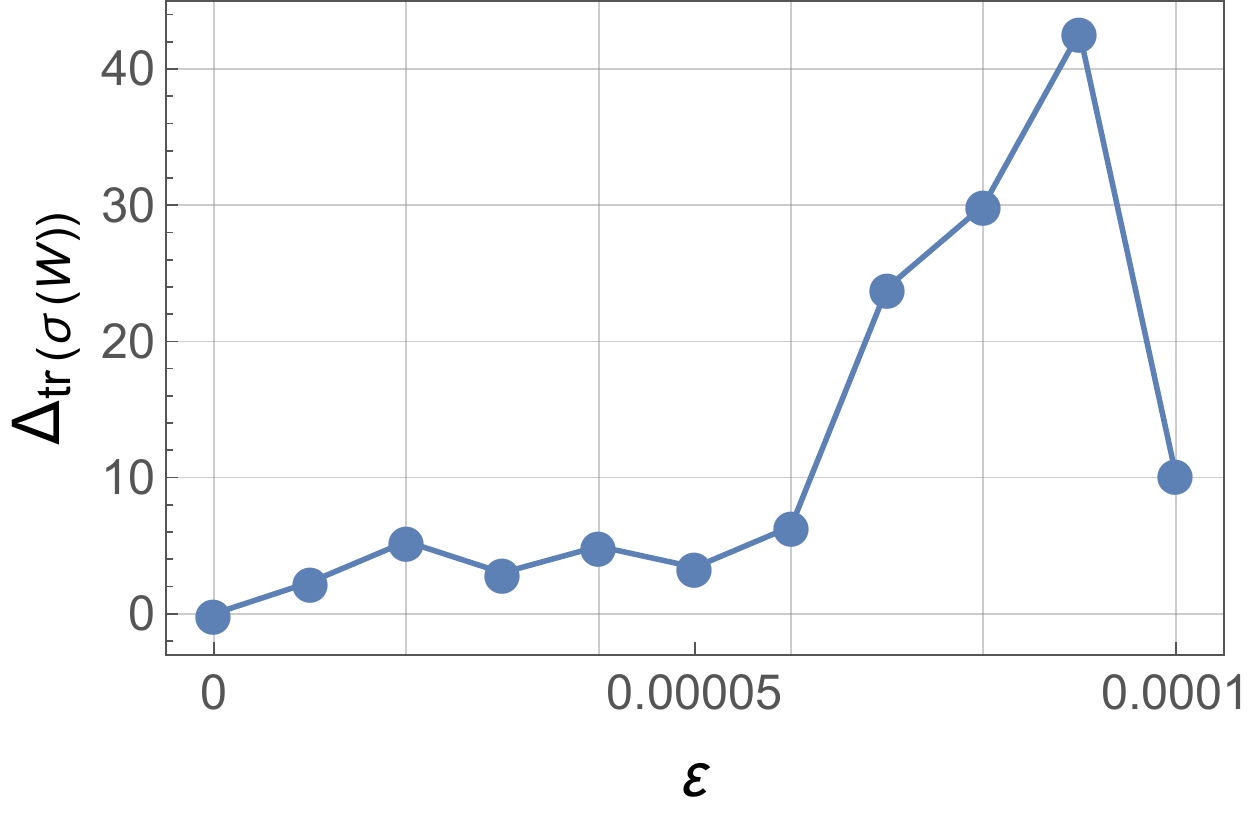}}
\caption{In~(a), a digraph with 10 vertices, $\mathcal G\equiv \mathcal G(A)$. Variation of the eigenvalues trace using the finite horizon Gramian matrix, with $k=10$, when we add uniformly generated noise with radius $\varepsilon$ only to the entry $A_{21}$. Notice that a small $\varepsilon$ value can reflect in a $\Delta_{tr(\sigma(W))}$ in the order of 100 000$\times \varepsilon$. }
\label{fig:exp_motiv_22}
\end{figure}

\vspace{2mm}
\setlength{\fboxsep}{0pt}
\fcolorbox{white}[HTML]{f3f3f3}{\parbox{.94\textwidth}{
\begin{remark}[Finite horizon Gramian]\label{rem:2}
Similarly, consider the use of the finite horizon controllability Gramian associated with $(A,B)$ and a finite horizon value $k>0$, given by
\[
 W=\sum_{i=0}^{k-1}A^i B B^\intercal (A^\intercal)^i,
\]
which enables to describe the controllability energy through its eigenvalues' sum. 
Additionally, consider the matrix $A$ as the adjacency matrix of the digraph depicted in Fig.~\ref{fig:exp_motiv_22}~(a), and \underline{add} uniformly generated random noise \underline{solely} to the entry $A_{21}=1$, where the noise is uniformly generated from the interval $]-\varepsilon,\varepsilon[$, for different $\varepsilon>0$ values ranging from $0$ up to $0.0001$. If $\varepsilon=0$ then we consider the absence of noise. 
In Fig.~\ref{fig:exp_motiv_22}~(b), for $k=10$, we illustrate the absolute value of the Gramian eigenvalues' sum variation for different values of $\varepsilon$, defined as $\Delta_{tr(\sigma(W))}=|tr(\sigma(W))-tr(\sigma(W_\varepsilon))|$, where $W_\varepsilon$ is the finite horizon Gramian with $k=10$ and matrix $A$ with noise $\varepsilon$ added to the entry $A_{21}$. 
Observe that a minimal perturbation to a single entry of matrix $A$ produces substantial changes in the Gramian eigenvalues' sum. 
Therefore, quantitative assessment of controllability energy might be misleading, and we should pursue other methods that enable us to assess controllability among other properties for large scale systems. 
\end{remark}
}}

Therefor, to avoid the types of problems illustrated in Remark~\ref{rem:1} and~\ref{rem:2}, we may adopt to study instead structural properties. 

In the context of structural systems, given the structural properties and assuming that the matrices' parameters are independent, we obtain  the following problem in terms of the generic controllability

\begin{equation}\label{eq:struct_1}
\begin{split}
    \mathop{\min}_{\Theta} & \quad f_\Theta(\bar A,\bar B)  \\
\text{s.t. }& grank(\mathcal C(\bar A,\bar B)) =n,
\end{split}
\end{equation}
where $grank$ is the \emph{generic rank} of $\mathcal C(\bar A,\bar B)$, i.e., the maximum rank that can be achieved with matrices $(A,B)$ that possess the same structural pattern as $(\bar A,\bar B)$.

The $grank$ notion finds a diversity of application, for instance, in Remark~\ref{Toeplitz}, we point a recent research line where its use has a central role in assessing the rank of Toeplitz matrices. 

\vspace{2mm}
\setlength{\fboxsep}{0pt}
\fcolorbox{white}[HTML]{f3f3f3}{\parbox{.94\textwidth}{
\begin{remark}[Toeplitz matrix]\label{Toeplitz}
Generic rank properties can also be used to assess the rank of matrices with predefined structure and parameter dependencies, as it is the case of Toeplitz matrices~\cite{reissig2013maximum}. An $n\times n$ matrix is a Toeplitz matrix whenever it can be defined as a matrix $A$ such that $A_{i,j} = c_{i-j}$, for constants $c_{1-n},\ldots, c_{n-1}$. 
For instance, for $n=3$ the following matrix $M$ is Toeplitz:
\[
M=\left[
\begin{array}{ccc}
    c_0 & c_{-1} & c_{-2}\\
    c_{1} & c_0 & c_{-1} \\
    c_2 & c_1 & c_0 \\
\end{array}
\right].
\] \hfill $\diamond$
\end{remark}
}}
\vspace{5mm}



We refer to the pair $(\bar A,\bar B)$ being structurally controllable if and only if $grank(\mathcal C(A,B))=n$ where $(A,B)$ have the structural pattern $(\bar A,\bar B)$. Thus, we can reformulate
~\eqref{eq:struct_1} as  the following structural optimization problem:
\begin{equation}\label{eq:struct_2}
\begin{split}
    \min & \quad g(\bar A,\bar B)  \\
\text{s.t. }& (\bar A,\bar B)\  \text{structurally controllable}.
\end{split}
\end{equation}
where $g:\{0,\star\}^{n\times n}\times \{0,\star\}^{n\times p}\rightarrow \mathbb{R}$. In the latter, it is easy to see that problem is combinatorial in nature given the domain of the objective function.

\subsubsection{Actuator placement} In the input design context, one often seeks to determine $\bar B $ given that $\bar A $ is assumed to be known and does not change. For instance, the objectives could be $\|\bar B\|_0$ that then models the problem of identifying the minimum number of state variables that need to be actuated to ensure structural controllability. Specifically, we seek to determine the solution to the following problem: given $\bar A$, find $\bar B$ that minimizes
\begin{equation}\label{eq:struct_3}
\begin{split}
    \min & \quad \|\bar B\|_0  \\
\text{s.t. }& (\bar A,\bar B)\  \text{structurally controllable}.
\end{split}
\end{equation}



\begin{figure}[!h]
\centering
\subfigure[]{\label{fig:a}\includegraphics[width=0.7\textwidth]{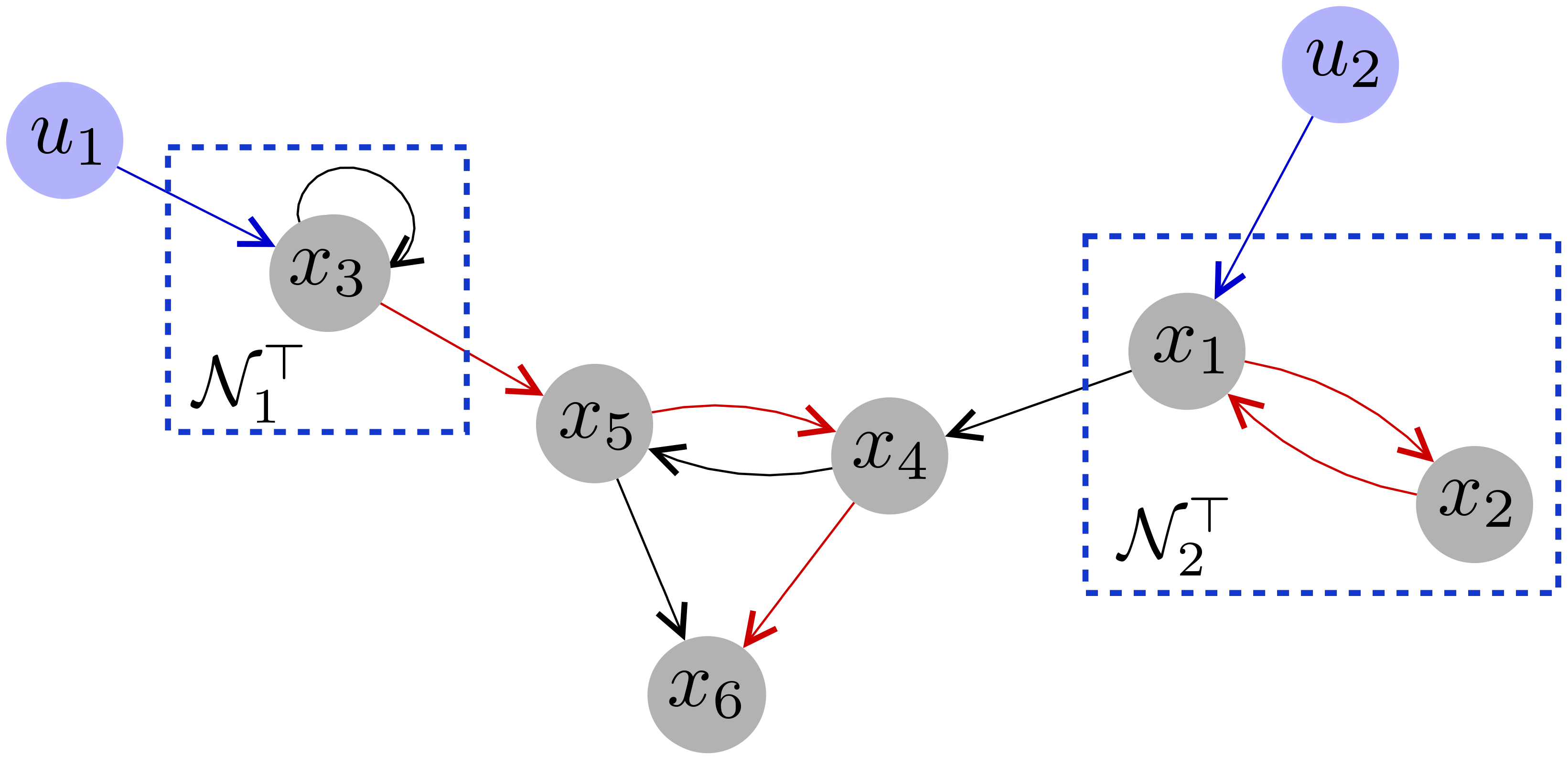}}
\subfigure[]{\label{fig:a}\includegraphics[width=0.63\textwidth]{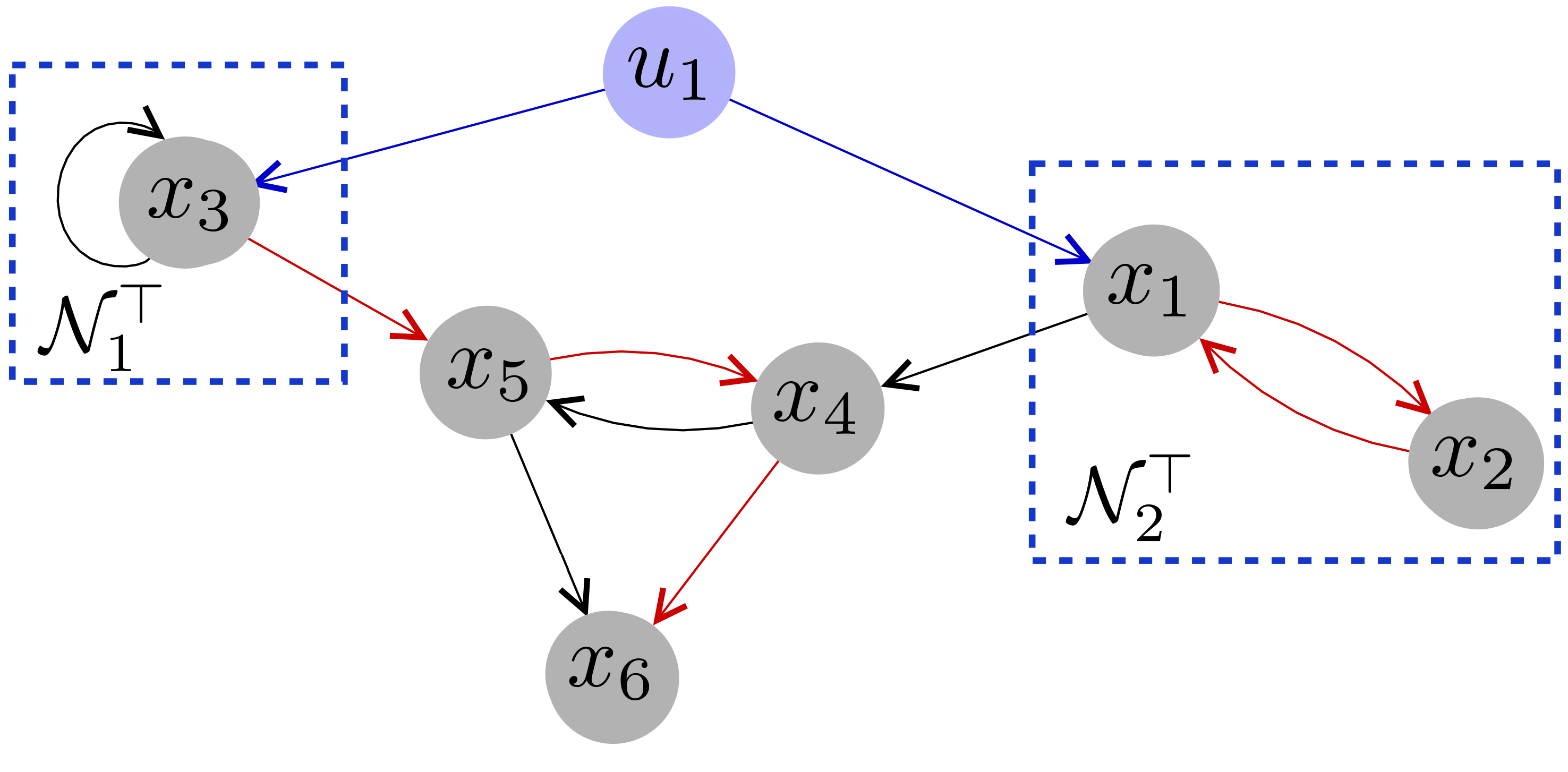}}
\caption{Digraph representation $\mathcal G(\bar A,\bar B)$ in the dedicated inputs scenario (a), and the non dedicated input scenario (b). The edges in red are the path and cycle decomposition obtained from the maximum matching of $\mathcal B(\bar A)$, depicted in Fig.~\ref{fig:bip}. The sets of vertices inside the dashed blue boxes are the source-SCCs and the vertex in the dashed green box is the target-SCC. Finally, the vertices in blue depict input variables. }
\label{fig:dedicated}
\end{figure}

To address the problem in~\eqref{eq:struct_3}, we can recast it to a combination of graph-theoretical procedures: maximum matchings on the state bipartite graph, and the unique decomposition of the state digraph into its \emph{strongly connected components} (SCC), i.e., disjoint subgraphs with the property that there exists a path between any two vertices in each subgraph. In particular, if we consider the additional constraint of at most one non-zero entry per column in $\bar B$ (i.e., when the inputs are dedicated), then the minimum number of dedicated inputs (i.e., non-zero columns) is given by~\cite{pequito2015framework}
\[
m = \max\{1,r+\beta-\alpha\},
\]
where $r$ is the number o right-unmatched vertices of the associated bipartite graph, $\beta$ is the number of source-SCC (i.e., connected components without incoming edges into their vertices, denoted by $\mathcal N^\top_1,\ldots,\mathcal N^\top_k$), and $\alpha$ is the maximum assignability index of the network. 
The maximum assignability index, $\alpha$, is the maximum number of source-SCC that contains right-unmatched vertices among all maximum matchings of the state bipartite graph. 
See the examples illustrated in Fig.~\ref{fig:dedicated} (a) in (b). 
Another key result, shown in~\cite{pequito2015framework}, is that, it is possible to compute a computationally efficient solution by reformulating the problem as a \emph{weighed} maximum matching problem~\cite{cormen2009introduction}.

Once a dedicated solution is determined (i.e., containing only dedicated inputs) all solutions to~\eqref{eq:struct_3} can be described as follows: one needs a nonzero entry in each of the rows corresponding to the state variables in right-unmatched vertices, whereas the remaining  $\beta$ nonzero entries can be arbitrarily distributed among the columns containing the nonzero entries in the previous step -- see illustrative example in Fig.~\ref{fig:dedicated} (b). Along the same lines, in~\cite{commault2015single}, the authors provide a polynomial complexity algorithm to achieve structural controllability for single-input systems.

Other actuator placement problems have been suggested and addressed in the literature. Among these problems, the oldest is that of determining the minimum number of inputs required to attain structural controllability, which equals the number $r$ of right-unmatched vertices of the state bipartite graph; thus, computationally efficiently solved; specifically, its time complexity is  $\mathcal O(n^3)$.

Other problems consider the minimum actuation cost on either the inputs~\cite{pequito2015minimum} or the state variables actuated~\cite{pequito2016minimum}; once again,  these problems can be solved by resorting to computationally efficient (i.e.,  $\mathcal O(n+m\sqrt{n})$, where $m$ denotes the number of non-zero entries of $A$)~\cite{olshevsky2015minimum}. In~\cite{doostmohammadian2019complexity}, the authors look into the particular case of self-damping systems.

Notwithstanding, such computational complexities might still be prohibitive when dealing with large-scale systems and, therefore,  in~\cite{faradonbeh2016optimality}, the authors proposed using fast maximum matching algorithms with optimality guarantees. Specifically, instead of incurring $\mathcal O (m\sqrt{n})$ time complexity to compute a maximum matching, it can be approximated with fast methods that hold linear time complexity ($\mathcal O(n)$).
Alternatively, random sampling schemes can be considered when the state space dimension is prohibitive even for polynomial complexity algorithms with certain approximation guarantees~\cite{ravandi2019controllability,jia2013control}. 
In~\cite{ravandi2019controllability}, the authors study the statistical characteristics of the actuator placement problem by randomly selection and assessing the resultant controllability properties of complex networks, whereas in~\cite{jia2013control}, a random sampling algorithm is developed to address the actuator placement problem. The algorithm not only provides a statistical estimate of the control capacity, but also to bridge the gap between multiple microscopic control configurations and macroscopic properties of the underlying network --
these results for the actuator placement problem complement some available heuristics to find a maximum matching for directed networks~\cite{chatterjee2013heuristic}.

\vspace{2mm}
\setlength{\fboxsep}{0pt}
\fcolorbox{white}[HTML]{f3f3f3}{\parbox{.94\textwidth}{
\begin{remark}
The (sparstest) minimum controllability problem~\eqref{eq:prob_min_control} given the matrix $A$ where the objective is to minimize the number of actuated variables, that is, $f_{\Theta}(A,B)=\|B\|_0$ is NP-hard~\cite{olshevsky2014minimal}. In contrast, the structural version of it in~\eqref{eq:struct_3} is polynomially solvable. In other words, from the definition of structural controllability, it follows that given an arbitrary matrix $A$ with the structure $A$ (assuming the parameters are independent of each other), the (sparstest) minimum controllability problem is almost surely polynomially solvable.
\end{remark}
}
}
\vspace{5mm}

\vspace{2mm}
\setlength{\fboxsep}{0pt}
\fcolorbox{white}[HTML]{f3f3f3}{\parbox{.94\textwidth}{
\begin{remark}
An important application of the actuator placement problem is that of \emph{leader selection} in the context of multi-agents~\cite{commault2013input, blackhall2010structural}. In this context, it is possible to assess the selection of leaders to ensure resilience as well as analytical properties that ensure the network to be resilient~\cite{jafari2011leader}. In addition, it is also possible to determine an efficient solution to the  \emph{distributed} leader selection that builds upon structural controllability tools, yet it guarantees (non-structural) controllability by determining proper weights in a fully distributed fashion~\cite{pequito2015distributed,tsiamis2017distributed}.
\end{remark}
}
}
\vspace{5mm}

Nonetheless, it would be wrong to assume that all structural controllability problems are equally difficult. Specifically, a small change in the assumptions of the optimization problem might lead to a NP-hard problem~\cite{pequito2015complexity}. For instance, consider the following optimization problem: given $\bar A$ and $\bar B$, determine the smallest subcollection of inputs $B(\mathcal I)$ ensuring structural controllability, where $\mathcal I=\{1,\ldots, p\}$ denotes the labels of the inputs described by the columns of the matrix $\bar B$, i.e.,
\begin{equation}\label{eq:struct_NPhard}
\begin{split}
    \min & \quad |\mathcal I|  \\
\text{s.t. }& (\bar A,\bar B(\mathcal I))\  \text{structurally controllable}.
\end{split}
\end{equation}

In fact, the decision version of the problem in~\eqref{eq:struct_NPhard} is NP-complete~\cite{assadi2015complexity}. Briefly speaking, this problem is as difficulty as several of the other well-known NP-complete problems. Nonetheless, it is possible to determine efficient approximations with provable optimality guarantees~\cite{moothedath2018flow} by  formulating a new graph-theoretic necessary and sufficient condition for checking structural controllability using flow-networks, and  proposing a polynomial reduction of the problem to the \emph{minimum-cost fixed-flow problem}, an NP-hard problem for which polynomial approximation algorithms exist.

\subsubsection{Sensor placement} This problem seeks to determine the sensors to be deployed, or sensing capabilities,
\begin{equation}\label{eqn:outputLTI}
    y(t)=Cx(t),
\end{equation}
with $y(t)\in\mathbb{R}^m$, 
required to attain \emph{structurally observability} -- a system is structurally observable if and only if there exists an observable $(A,C)$ with the structural pattern $(\bar A,\bar C)$. By invoking the duality between controllability and observability in the context of LTI systems, it readily follows that all the actuator placement problems can be posed as sensor placement problems. 

Notwithstanding, a variety of sensor placement problems have been proposed in~\cite{commault2007sensor,boukhobza2007state}.  In~\cite{liu2018partial}, the authors proposed to determine a collection of sensors to ensure that a certain observability subspace allow to retrieve certain state variables of interest.    Additional, sensor placement can be achieved under possible cost constraints~\cite{liu2017scheduling,liu2019configuration}. 
 
The authors, in~\cite{doostmohammadian2017observational}, proposed to leverage the structure of the topology by performing contractions of the graph (i.e., to obtain topological equivalent classes). By doing this, they improved the efficiency of the  sensor deployment algorithm.

In~\cite{commault2011sensorfeedback}, sensor placement is considered towards disturbance rejection by measurement feedback.

Additionally, sensor placement can be considered to attain resilience/robustness with respect to sensor failure that maintain structural observatibility~\cite{commault2008observability,boukhobza2009state,boukhobza2010partial}. 
Lastly, in~\cite{kruzick2017structurally}, the authors address the observability problem when backbone nodes (e.g., routers) are considered in the context of sensor networks.

\subsubsection{Dynamics topology design problem to attain structural controllability} In contrast with the actuator placement problem where the goal is to add actuation capabilities, we now focus on changing the dynamics topology. Specifically, we have the following problem: given $(\bar A, \bar B)$, find the sparsest structural perturbation

\begin{equation}\label{eq:struct_dyn}
\begin{split}
    \min & \quad \|\bar \Delta\|_0  \\
\text{s.t. }& (\bar A+\Delta,\bar B)\  \text{structurally controllable}.
\end{split}
\end{equation}

The problem in~\eqref{eq:struct_dyn} can be efficiently solved~\cite{chen2018minimal,zhang2019structural,mu2019guaranteed}. Alternatively, we can consider a combination of both actuator placement and dynamics design, as well as possible costs associated with these, as explored in~\cite{zhang2019minimal}. 
 Given a structured system, in~\cite{zhang2019minimalPerturbations} it is proven that  it is NP-hard to add the minimal cost of links,
including links among state variables (i.e., state links) and links from the existing inputs to state variables (i.e., input links), from a given set of links to make
the system structurally controllable. Similarly, it is NP-hard to determine the minimal cost of links whose deletions deteriorate structural controllability of the
system, even when the removable links are restricted in either the input links
or the state links. In other words, we can assess how many links to remove to withhold the controllability properties~\cite{zhang2019minimal,dey2018minimum}. 
The latter problem is associated with the resilience and security of structural systems, where the problem is that of determine the impact of edges failure/removal in the structural controllability~\cite{ramos2015analysis,jafari2010structural,rahimian2013structural,alcaraz2014recovery}. 

A natural extension of the problem~\eqref{eq:struct_dyn} is that of considering the dynamics topology design in the context of networked dynamical systems. Specifically, consider a set of $N$ subsystems $\{\bar A_i,\bar B_i\}_{i=1}^N$ that are interconnected through $\bar A_{i,j}$ with possible restrictions on its own as they represent outputs from subsystem $j$ that serve as inputs in the subsystem $i$. 
In this context, if all the systems are alike, then we have \emph{homogeneous} networked dynamical systems, otherwise, we have \emph{heterogeneous} networked dynamical systems. In both cases, as these systems increase their size, it is of interest to determine easy to verify conditions for structural controllability that depend only on the interconnection structure and possibly on subsystems general properties~\cite{carvalho2017composability}. In fact, due to the geographically distributed nature of these system, it is often desirable that structural controllability is assessed in a distributed fashion~\cite{carvalho2017composability}.

Subsequently, we can pose the problem of determining $\Delta_{i,j}$ (with the same dimensions of $A_{i,j}$) under possible structural restrictions that will change the interconnection between the different subsystems. Unfortunately, some of the problems are NP-hard, which than require the development of approaches that will efficiently provide an approximate solution with probably optimality guarantees~\cite{moothedath2019optimal}.

\subsubsection{Joint dynamics and input design to attain sctructural controllability}

Given a structured system, in~\cite{zhang2019minimalPerturbations} it is proven that it is NP-hard to add the minimal cost of links,
including links among state variables (i.e., state links) and links from the existing inputs to state variables (i.e., input links), from a given set of links to make
the system structurally controllable, and similarly, it is NP-hard to determine the minimal cost of links whose deletions deteriorate structural controllability of the
system, even when the removable links are restricted in either the input links
or the state links.

\subsection{Dynamics Topology Design:  Structural Stabilizability and Structural Decentralized Control}\label{sub:dyn_top_des}

In Section~\ref{sec:systemDesignsC}, we focus in optimization problems that aim to attain structural controllability. In this section, we focus on structural stabilizability and decentralized control.

\subsubsection{Structural Stabilizability}\label{sec:strucStab} A system is said to be structurally stable if the dynamics described by $A$ is stable (i.e., its eigenvalues lie within the unit circle in the complex plane) for almost all matrices that have a given structural pattern $\bar A$. This problem as been addressed in~\cite{belabbas2013sparse,kirkoryan2014decentralized}, where the authors study the patterns of the matrices that are stable.

These  can be related with the control configuration design problems presented next.

\subsubsection{Control configuration design problem (or, structural feedback selection)}

In the context of decentralized control, it often happens that the outputs of the system~\eqref{eqn:outputLTI} are not available for feedback to all the actuators of the system. Specifically, when we consider \emph{static output feedback}, where the control law takes the form
\begin{equation}
    u(t)=Ky(t)
\end{equation}
where the \emph{gain}  $K\in\mathbb{R}^{p\times m}$ is time-invariant. In particular, the gain satisfies an \emph{information (structural) pattern} $\bar K$, where an entry $\bar K_{ij}=\star$ if the output (or sensor) $j$ is available to the input (or controller) $i$. In this context, it suffices one entry $\bar K_{ij}=0$ to be in the scenario of decentralized control. Notice that this generalizes the initial notion of decentralized control where the information pattern is taken to be diagonal or block-diagonal. As such, when the information pattern is a full matrix we are in the centralized control setting.

When designing closed-loop control laws, we often seek to shape the performance of the system through the so-called pole placement, that is we seek to design the gain such that the spectrum of the closed-loop systems given by $A+BKC$ is as specified. When we considered \emph{state feedback} (or alternatively, when the output matrix equals the identity) it is possible to attain arbitrary pole place when the system is controllable. 
Nonetheless, when we consider output feedback arbitrary pole placement may not be possible to attain. In fact, in the context of decentralized control, we often seek to guarantee that the system does not have (decentralized) \emph{fixed modes}, that result in eigenvalues in the spectrum of the closed loop system that do not change regardless of the gain chosen that satisfies a specified information pattern. Specifically, let $\mathcal K=\{K:K_{ij}=0 \text{ if } \bar K_{ij}=0\}$, then the system described by $(A,B,C)$ has fixed modes with respect to the information pattern $\bar K$ if and only if
\[
\bigcap\limits_{K\in\mathcal K} \sigma(A+BKC) \neq \emptyset,
\]

\noindent where $\sigma(M)$ denotes the spectrum of the square matrix M. It is also possible to consider fixed modes in the context of structural systems, which are referred to as \emph{structurally fixed modes}, and are defined as follows~\cite{sezer1981structurally}: a system $(\bar A,\bar B,\bar C)$ has structurally fixed modes with respect to the information pattern $\bar K$ if there exists a realization $(A,B,C)$ with the same structural pattern of $(\bar A,\bar B,\bar C)$ that has fixed modes with respect to the information pattern $\bar K$. In other words, the system has no structurally fixed modes if there is no realization satisfying the structural pattern of the state space matrices that yields fixed modes with respect to a given information pattern. In particular, it is easy to see that if a system is not controllable, then the uncontrollable modes are also fixed modes – and similar consequences can be drawn in the context of structural systems. Yet, controllability (structural controllability, respectively) do not guarantee the non existence of fixed modes (structurally fixed modes, respectively) with respect to a specific information pattern. In fact, in~\cite{lee2017structurally}, the authors show that structural fixes modes cannot co-exist with quadratically invariant or partially nested information structures~\cite{mahajan2012information}.

That said, it is often the case that one seeks to find the solution to the following optimization problem: given the system’s structure $(\bar A,\bar B,\bar C)$, find the sparsest information pattern $\bar K$ as follows
\begin{equation}\label{eq:struct_cc}
\begin{split}
    \min & \quad \|\bar K\|_0  \\
\text{s.t. }& (\bar A,\bar B,\bar C,\bar K)\  \text{has no structurally fixed modes}.
\end{split}
\end{equation}

The interest on the problems  of the form~\eqref{eq:struct_cc} go back to~\cite{trave1987minimal}, where possible different costs could be considered. Nonetheless, it was only recently that their complexity, as well as efficient algorithms to compute the solutions or approximations in the case where the problem is NP-had, were addressed in~\cite{pequito2015framework,pequito2015minimum,moothedath2018minimum}. The authors, in~\cite{moothedath2019sparsest}, propose polynomial-time algorithms to design the sparsest feedback gain for cyclic systems, and  in~\cite{pequito2018analysis}, efficient algorithms for the sparsest robust design for cyclic systems are described. 

It is worth to mention that there are two natural extension to the above problems: (i) in practice, it could be possible that one seeks to guarantee that the system only has stable fixed modes. In other words, we seek to determine structural feedback links that ensure that the system has no unstable structural fixed modes; and (ii) instead of a static (memoryless) feedback controller, we can equip the controller with memory and in that case it might be possible to ascertain structural stability without the need to guarantee that the system $A+BKC$ is structurally stable -- see Section
~\ref{sec:strucStab}. Towards this direction, the authors in~\cite{pajic2011wireless,pajic2013topological}, leverage the sensor network that has it own dynamics and, subsequently, behaves as a controller with memory, which in closed-loop is able to attain stability conditions of the composed system -- see also Section~\ref{sub:wir_net} for the application in wireless sensor networks.
Lastly, in~\cite{kalaimani2013generic}, the authors explored the minimal controller structure for generic pole placement.

\vspace{2mm}
\setlength{\fboxsep}{0pt}
\fcolorbox{white}[HTML]{f3f3f3}{\parbox{.94\textwidth}{
\begin{remark}[Connection between structural and non-structural]\label{structNonStruc}
It is also possible to connect some of the above ideas originated in structural systems theory with non-structural results. Specifically, on the connection between the structure and the results for a parametrized system, in~\cite{torres2015graph}, the authors propose  conditions of single input-output feedback stabilization. In~\cite{lee2017decentralized}, the authors provide a simple design mechanism to obtain a gain to attain pole placement with arbitrary structure of a dynamic output feedback.
\end{remark}
}
}
\vspace{5mm}

\subsubsection{Co-design of actuator/sensor placement and control configuration}

Due to the inter-dependency between the actuator/sensor placement problem and the control configuration problem, a more interesting problem that several authors have addressed  is that of the \emph{co-design} problem of actuator/sensor placement and control configuration, where the goal is as follows: given $\bar A$ determine $(\bar B,\bar C,\bar K)$ such that
\begin{equation}\label{eq:struct_coDesign}
\begin{split}
    \min & \quad \|\bar B\|_0+\|\bar C\|_0+\|\bar K\|_0  \\
\text{s.t. }& (\bar A,\bar B,\bar C,\bar K)\  \text{has no structurally fixed modes}.
\end{split}
\end{equation}

The problem in~\eqref{eq:struct_coDesign} was addressed in~\cite{pequito2015framework}, which can be polynomially solvable. 
Extensions that consider heterogeneous costs on both the actuator/sensor placement and the control configuration lead to NP-hard problems, but some subclasses of systems can still be solvable polynomially~\cite{pequito2015minimum}. 
Remark that  if a problem is NP-hard it does not mean that all instances of such problems are equally difficult to solve nor to approximate. Similarly to the actuator placements, if we seek to determine the sparsest information pattern that uses the smallest collection of inputs/outputs from a specified collection of possible inputs/outputs, then the problem is also NP-hard, and subsequently, several  schemes should be considered to efficiently obtain a suboptimal solution with possibly optimality guarantees~\cite{moothedath2018minimum}.

\newpage

\subsection{Controllability Index}\label{{subsub:cont_index}}

The controllability index plays a major role in assessing the minimum number of time steps that lead to the existence of a control law capable of steering the system state towards a desired goal. Specifically, the controllability index $k^*$ is given by
\begin{equation}
k^*=\arg \min  \{{k\in\mathbb{N}}: rank \ \mathcal C_k(A,B)= n\}, 
\end{equation}
where $C_k(A,B)=[B\; AB\; A^2B\;\ldots\;A^{k-1}B]$ is the partial controllability matrix.

In~\cite{pequito2017trade}, the authors provide a characterization of the structural controllability index, for which they provide a near-optimal approximation scheme to discover the minimum number of inputs required to attain a given controllability index.  In~\cite{ding2016optimizing}, the authors propose to minimize the controllability index upon a budget on the number of inputs. More recently, in~\cite{ramos2020generating}, the authors proposed a  generative model to attain a specific \emph{actuation spectrum}, which captures the trade-offs between the minimum number of state variables required to attain structural controllability in a given number of time-steps.

\subsection{Target Controllability (Output Controllability)}\label{subsub:target_control}

Consider a system $(A,B)$ and a target node set $\mathcal C=\{c_1,\hdots, c_S\}$. Let $C(\mathcal C)$ denote the $n\times n$ matrix that comprises the $S$ rows of the $n\times n$ identity matrix indexed by $\mathcal C$. 
The triple $(A,B,C)$ is target controllable if 
\[rank([\begin{array}{cccc}CB & CBA & \ldots & CBA^{n-1}\end{array}])=S.
\]

The seminal work of Gao \emph{et al.}~\cite{gao2014target} introduced the target control problem and proposed a $k$-walk theory to address it for directed-tree like networks with a single input.

In~\cite{van2017distance}, using the structural output controllability properties, the authors study target controllability of dynamic networks~\cite{monshizadeh2015strong}. 
In~\cite{czeizler2018structural}, the authors provide more efficient algorithms than the ones introduced in~\cite{gao2014target} and illustrate some of the limitations of the $k$-walk approach. 

In~\cite{li2019target}, the authors present a method to allocate a minimum number of external control sources that assure the target controllability of LTI systems complemented with linear observers. The method locates a set of directed paths and cycles to cover the target set. 

In~\cite{moothedath2019target}, the target controllability of a network, i.e., a family of structured systems, is investigated. The authors displayed a bipartite matching-based necessary condition to discover a network's target controllability. 
The decision of which state variables should be controlled such that only a subset of state variables is controllable was proved to be, in general,  NP-hard~\cite{gao2014target}. 
Additionally, if the selection of variables is constrained to a pre-defined set, the problem is also NP-hard~\cite{guo2017constrained}. 
Remarkably, the former problem is polynomially solvable for symmetric state matrices~\cite{li2018structural}. 
The computational complexity of the target structural controllability problem is studied in~\cite{czeizler2018fixed}. The authors show that the problem admits a polynomial-time complexity algorithm when parameterized by the number of target nodes. In general, they show that it is hard to approximate at a factor better than $\mathcal O(\log n)$. 
The work in~\cite{commault2019functional} considers a stronger controllability notion, when only a small set of state variables can be actuated. 
Instead of using the Kalman controllability, the authors require the ability to drive the target variables as time functions. 
In this setup, they solve the described controllability problem using a functional approach. 
Lastly, in~\cite{guan2019target}, the authors propose sufficient conditions for target controllability under possible switches in the topology. 

\subsection{Edge dynamics}\label{subsub:edge_dyn}

A different research direction considers, instead of controlling the nodes, controlling the edges of the graph. In this scenario, notions of controllability are possible using the notion of dual graph~\cite{nepusz2012controlling,pequito2016structural,shen2018structural}.

\subsection{Structural identifiability}\label{subsub:Struct_ident}

Consider a LTI system whose state space matrices are parametrized as $(A(\theta),B(\theta),C(\theta))$, and the corresponding transfer function $G(s,\theta)=C(\theta)(sI-A(\theta))B(\theta)$, then structural identifiability can be defined as follows~\cite{van2008determining}:
A model structure $G(s,\theta)$ is (locally) structurally identifiable at the optimal set of parameters $\theta^*\in\Theta$, where $\Theta$ is the set of all possible parameters, if for all $\theta_1,\theta_2\in\Theta$ in the neighborhood of $\theta^*$ and for all $s$ one has
\[
G(s,\theta_1)=G(s,\theta_2)\implies \theta_1=\theta_2.
\]
Briefly speaking, an identifiable model structure yields distinct models with distinct parameters, and if distinct parameters yield the same model (i.e., $G(s,\theta_1)=G(s,\theta_2)$ when $\theta_1\neq \theta_2$), then the model structure is not identifiable.

A diverse set of conditions for the structural identifiability have been proposed in~\cite{agbi2012parameter,agbi2014decentralized,canto2009structural}.   
In~\cite{canto2011identifiability} the authors addressed the identification of discretization schemes for partial differential equations. 
In~\cite{bhela2017enhancing,lagonotte1989structural}, the authors assess the identifiability of a linearized power grid system model, by leveraging generic properties. Lastly, in~\cite{ideta2004structural}, the authors presented conditions  for the structural identifiability of singular systems with delay.


\subsection{Input–output decoupling structure}\label{subsub:IO_decouple_struct}

The problem of input–output decoupling is also known as row-by-row decoupling problem or Morgan's problem. 
Given the system in~\eqref{eq:lti}, consider that we have the same number of input and output (i.e.,  $B$ and $C$ have the same dimension). 
The goal of this problem is to find a state feedback of the form $u(t)=F x(t)+J v(t)$, with $J$ nonsingular, such that the closed-loop system transfer matrix
\[
T_{F,J}(s) = (C+DF)(sI-A-BF)^{-1}BJ+DJ
\]
is diagonal and nonsingular. 

The work in~\cite{abad2014graph} finds a relationship between the network's graph topology and the infinite-zero and finite-zero structures of an input-output network's dynamics. The authors express the zero-dynamics state matrix as a perturbation of the reduced graph matrix (a sub-matrix of the network dynamics state matrix). 
In~\cite{conte2019invariance}, an original approach to studying a class of structured systems was introduced. The authors propose the following new graph-theoretic notions: invariance, controlled invariance, conditioned invariance, and essential feedback. 







\subsection{Fault detection and isolation (FDI)}\label{subsub:DFI}

The fault detection and isolation problem (FDI) consists of designing a set of signals, for instance,  via observers, called residuals. 
This set of signals should be such that the transfer matrix from the disturbance and control inputs to the residuals is zero. 
Moreover, the transfer matrix from the faults to the residuals must have a specific form, e.g., it should be diagonal or triangular. 
These residuals are insensitive to controls and disturbances but sensitive to faults. This property allows them to detect and isolate the faults. 

In~\cite{commault2002observer}, the authors provide necessary and sufficient conditions under which the FDI problem has a solution for almost any values of the free parameters. These conditions are expressed in terms of input-output paths of the directed graph associated with the original LTI system (not the closed-loop LTI system associated with the residual error dynamical system). 
Alternatively, in~\cite{commault2008structural}, the authors use structural systems to analyze the sensor location in the Fault Detection and Isolation problem. 
This analysis allows determining the minimal number of required extra sensors and the sets of internal variables that need to be measured for solving the FDI problem. 
Decentralized FDI is studied in~\cite{sauter2006decentralized}, and the authors use structural control to draw necessary and sufficient conditions for detectability and isolation. 
In~\cite{commault2011sensor,chamseddine2009optimal}, the authors addressed the problem of sensor placement to ensure that FDI can be performed, when disturbances may exist.  
Additionally, in~\cite{boukhobza2008graph}, the authors study the FDI problem but for the case of bilinear systems.  
For an overview of the different approaches to tackle the DFI problem, we refer the reader to the work in~\cite{simon2013reliability}. 
Lastly, in~\cite{staroswiecki2007structural}, the authors present a structural view of FDI. The authors show that, by resorting to structural analysis, it is possible to establish a link between critical faults and reliability. 



\subsection{Security and Resilience}\label{subsub:security}



Distributed control systems (DCS) usually rely on different components. Therefore, they are exposed to malicious attacks.  
The use of DCS in critical infrastructures makes their security an important issue. 
Nefarious incidents connected to the security of DCS include the Stuxnet uranium plant attack~\cite{langner2011stuxnet} and the Maroochy Shire~\cite{abrams2008malicious} episode. 
Consequently, there has been a growing effort to mitigate  DCS from being exposed to undetectable attacks. 
In~\cite{sundaram2008distributed}, the authors design a scheme that allows nodes of time-invariant connected networks to attain consensus on any arbitrary function of the initial nodes' values in a finite number of steps for almost all weight matrices with the same structure. 
In~\cite{liu2011false}, for the power grid's design, the authors render algebraic conditions that allow an adversary to generate state estimation errors.  
The work in~\cite{sandberg2010security} suggests multiple security indices for sensors, allowing a system operator to identify sparse power grid attacks.  
In~\cite{mo2010false}, the bias that an undetectable adversary may introduce into the state estimation error of control systems and sensor networks is studied. 
Afterward,~\cite{sundaram2010distributed} and~\cite{pasqualetti2011consensus} address the resilience of consensus-based algorithms. 
The first work determines graphical conditions under which a set of agents can compute a function of their initial states in the presence of malicious nodes. 
The second work, using connectivity and left invertibility, delineates attack identifiability and detectability. 
Subsequently, in~\cite{sundaram2010wireless}, the authors study the design of an intrusion detection scheme for DCS, which identifies malicious agents and can recover from the attacks.

In the context of structural observability, we also have to keep in mind that the properties hold generically. 
Consequently, we can explore the set with zero Lebesgue measure to design attacks that will not be identifiable from the observer perspective~\cite{pasqualetti2013attack}. 
In other words, structural observability is only a necessary condition. In particular, the previous authors explore the left-invertibility to regular descriptor systems. 
In the same line, in~\cite{weerakkody2016graph}, it is considered  the secure design problem in the context of distributed control systems to guarantee the detection of stealthy integrity attacks. 
In~\cite{milovsevic2018security}, the authors propose a security index upon the previous definitions. 
Also, in the same research direction, in~\cite{weerakkody2017robust}, the authors propose to explore structural observability properties to prevent zero dynamics attacks. 
In~\cite{milovsevic2020actuator}, the authors leverage structural systems theory to derive an upper bound and the robust index.





The assess of robust control in the context of resilient leader selection, and analytical properties that ensure the network is resilient is performed in~\cite{jafari2011leader}.

In~\cite{zhang2019driver}, the authors propose to evaluate the minimum number of additional actuation capabilities needed to ensure structural controllability under possible failures in the form of a security index. In~\cite{alcaraz2017resilient,alcaraz2017cyber}, the authors explore self-healing properties to guarantee structural controllability through various centrality measures and severity degrees. 

The sensor placement problem can be considered to attain resilience/robustness with respect to sensor failure that maintain structural observatibility~\cite{commault2008observability,boukhobza2009state,boukhobza2010partial}. 

In~\cite{pequito2018analysis}, the authors provided efficient algorithms for the sparsest robust feedback design for cyclic systems. 
In like manner, to assess security properties, several criteria that serve as a denial of service were deemed  in~\cite{ramasubramanian2016structural}. 

Additionally, in~\cite{jafari2010structural}, the authors assess the impact of link failures into structural controllability (i.e., zeroing the free parameters of the autonomous system matrix). In~\cite{maza2012impact,dakil2015disturbance}, the authors address the issue of ensuring structural controllability under the scenario where the actuators can fail with known probabilities. Similarly, in~\cite{liu2016sensor,dakil2015generic}, the authors propose a methodology to determine the sensors which ensure structural observability within a certain probability. 
In~\cite{guan2019target}, the authors present sufficient conditions for target controllability under possible switches in the topology. 

In~\cite{zhang2018efficient,zhang2017iterative}, the authors assess the impact of node removal on structural controllability. They propose a strategy to recalculate the minimum number of inputs to guarantee structural controllability. 
In~\cite{shoukry2015imhotep}, the authors introduce the notion of structural abstraction to address formal guarantees in the context of the security index of systems.

In the context of hybrid systems, the authors of~\cite{ramos2015analysis} introduced a tool for the design and verification of structural controllability of theses systems, and they utilized this tool in the context of the 
analysis and design of electric power grids, ensuring robustness to the link failures~\cite{ramos2015analysis}. 
The authors in~\cite{jafari2010structural,rahimian2013structural,alcaraz2014recovery} also consider the impact of edge failure/removal to ensure resilient structural controllability.

\subsection{Privacy}\label{subsub:privacy}

We start by noticing that privacy in the context of dynamical systems is intrinsically intertwined with the notion of observability. Specifically, states are referred to as private if they do not belong to the observability subspace. Subsequently, by increasing sensing capabilities is more likely to increase the dimension of the observability subspace and the systems topology constraints the latter. Thus, reducing the privacy of the system by increasing the number of state variables possible to recover from the systems' outputs.

In~\cite{yan2012distributed}, the authors explore the concept of observability as a way to retrieve the objective function used by different agents in a network that leverages the knowledge of the different iterations performed by an iterative algorithm. In~\cite{pequito2014design}, the authors propose to design the network topology to ensure that some nodes' state cannot be retrieved by some of the agents in the network. 
The central idea is that a star network leads to a scenario where the center node can retrieve the states of all the corner nodes, and, in turn, the corner nodes can retrieve the center node's state, but not the remaining corner nodes' states. 
Intuitively, we should create a network contemplating a trade-off between the number of start networks and the connections between these networks.  
Some of these insights were further considered in the context of wireless sensor networks~\cite{pequito2015smart}. In a different direction, in~\cite{lin2015decentralized}, the authors consider privacy-preserving decentralized matrix completion, in which a network of agents collaborate to complete a low-rank matrix that is the collection of multiple local data matrices.

\section{Other Subclasses of Structural Systems}\label{sec:other_sublasses}

In what follows, we present an overview of recent results in the context of super classes of linear time-invariant systems, i.e., classes that include LTI systems, to which several structural conditions have been proposed, as well as the corresponding design problems.

\subsection{Composite linear-time invariant systems}\label{sec:compositeSystems}

Consider $r$ continuous-time LTI systems described as
\[
\dot x_i(t) = A_i x_i(t) + B_i u_i(t),\,i=1,\ldots,r,
\]
where the state $x_i(t)\in\mathbb R^{n_i}$ and the input $u_i(t)\in\mathbb R^{p_i}$. 
By considering the interconnection between subsystem $i$ and $j$, for all possible subsystems, we collect the interconnected dynamical system represented as follows:
\[
\dot x(t) = \underbrace{\left[\begin{array}{cccc}
    A_1 & E_{1,2} & \cdots & E_{1,r} \\
    E_{2,1} & A_2 & \cdots & E_{2,r} \\
    \vdots & \ddots & \ddots & \vdots \\
    E_{r,1} & \cdots & E_{r,r-1} & A_r
\end{array}\right]}_{A}x(t)+\underbrace{\left[\begin{array}{cccc}
    B_1 &0 & \cdots & 0 \\
    0 & B_2 & \cdots & 0 \\
    \vdots & \ddots & \ddots & \vdots \\
    0 & \cdots & 0 & B_r
\end{array}\right]}_{B}u(t),
\]
where the state is $x(t)=[x_1^\intercal(t)l\,\cdots\,x_r^\intercal(t)]\in\mathbb R^{n}$, $n=\sum_{i=1}^r n_i$, and $u(t)=[u_1^\intercal(t)\,\cdots\,u_r^\intercal(t)]\in\mathbb R^{p}$, $p=\sum_{i=1}^r p_i$. 
Additionally, $E_{i,j}\in \mathbb R^{n_i\times n_j}$ is the connection matrix from the $j$th to the $i$th subsystems. 

Although these are linear systems, due to their wide applicability we consider it as a class on its own as the tools and characterizations differ from those used in general linear time-invariant systems. For instance, in~\cite{carvalho2017composability,xue2019structural}, the authors leverage the similarity of the graph to derive easy to verify necessary and sufficient conditions for structural controllability for homogeneous systems (i.e., similar structure of subsystems) and serial systems to be structurally controllable. 
They design distributed algorithms to verify necessary and sufficient conditions to assure structural controllability for any interconnected dynamical system, consisting of LTI subsystems, in a distributed fashion. 
In~\cite{xue2019structural}, it is shown that structural controllability for networks involving homogeneous subsystems may not decompose into subsystem-level and network-level conditions since the system can have structural network invariant modes. 
Similarly, in~\cite{doostmohammadian2018structural}, the authors analyze the computational cost of sensor networks optimization monitoring structurally full-rank systems under distributed observability constraints. 
In~\cite{moothedath2019optimal}, the authors introduced methods to identify a minimum cardinality set of interconnection edges that the subsystems of a heterogeneous system should ascertain between them to yield a structurally controllable composite system.

In~\cite{commault2018classification,commault2019structural,doostmohammadian2019minimal,wang2017structural}, the authors render a classification of subsystem nodes based on their role in the overall structural controllability.  
Moreover, we can assess the minimum number of dedicated inputs required to accomplish structural controllability properties for subclasses of systems such as bipartite networks~\cite{nacher2013structural}, and in multiplex networks~\cite{nacher2019finding}. 
More recently, in~\cite{bai2019block}, the authors provide a divide and conquer strategy that divides the system into different blocks to address the minimum input design problem for structural controllability for large-scale systems. In~\cite{moothedath2020optimal}, the authors explore the problem of designing composite systems by introducing two indices (i.e., maximum commonality index and dilation index) that explore the trade-offs on the size of subsystems and connections among these.

\subsection{Descriptor linear time-invariant systems}\label{subsub:Descriptor_LTI}

Descriptor linear time-invariant systems (also referred to as implicit linear time-invariant systems) are those of the form:
\[
Ex(t+1)=Ax(t)+Bu(t), \ t=0,1,\ldots,
\]
where the matrix $E$ has appropriate dimensions, and the structural description of the systems is given by $(\bar E,\bar A,\bar B)$.

Although they represent a small modification concerning the linear time-invariant systems, their descriptive power is way broader. The analytical tools to assess these systems properties (e.g., controllability) are somewhat different. 
Nevertheless, some cases are more straightforward and resemble those of linear time-invariant systems, as it is the case of regular descriptor systems~\cite{lewis1992tutorial}. Subsequently, all the problems overviewed in the context of linear time-invariant systems can be posed in the context of descriptor time-invariant systems. Remarkably, only a few papers address some of these during the period cover in this overview. 
For instance, in~\cite{boukhobza2006observability}, the authors provided conditions for structural observability of descriptor systems. In~\cite{clark2017input}, the authors addressed the actuator placement problem to ensure structural controllability. 
Lastly, in~\cite{mathur2018design}, the closed-loop properties are explored for low dimensional descriptor systems.

\subsection{Linear time-invariant systems with  delays}\label{subsub:delay_LTI}

In this case, the main focus has been in obtaining conditions of structural controllability when dealing with linear time-invariant systems with delays~\cite{qi2016structural,van2018structural}.

\subsection{Linear time-invariant systems with  unknown inputs}\label{subsub:unknownInput_LTI}

Another class of linear time-invariant systems for which structural systems properties have been explored is that with unknown inputs~\cite{commault2001unknown,boukhobza2014graph}. Specifically, in~\cite{commault2001unknown}, the authors explored the unknown input observers design analysis from a structural systems perspective. 
In~\cite{boukhobza2014graph}, necessary and sufficient conditions which guarantee that a given set of unknown parameters describing the system's model is structurally identifiable are drawn in graphical terms.

\subsection{Bilinear systems}\label{sub:bilinear}

A \emph{bilinear system} can be formally described as~\cite{mohler1980overview,isidori1974realization}
\[
x(t+1)=Ax(t)+Nx(t)u(t)+Bu(t), \ t=0,1,\ldots,
\]
where $N$ is a matrix with appropriate dimensions and $u(t)$ is a scalar input.

Structural controllability properties are addressed in the context of homogeneous (no linear component of the input, i.e., $B=0$) rank-1 bilinear systems~\cite{ghosh2016structural,ghosh2017graphical}. In~\cite{tsopelakos2018classification}, the authors addressed the structural controllability of driftless bilinear control systems. Additionally, they studied the accessibility of these with a drift.

In~\cite{boukhobza2007observability,boukhobza2008generic,canitrot2008observability}, the authors scrutinize the structural observability of bilinear systems, and in~\cite{boukhobza2008graph}, the authors proposed the study of sensor selection in the context of FDI. 
Lastly, in~\cite{svaricek2006discussion}, the author presents an interesting discussion on the use of structural systems to assess uniform observability of bilinear systems.

\subsection{Discrete-time fractional-order Systems}\label{sub:frac_order}

Fractional order systems are successfully used to model different physiological processes (e.g., trains of spikes, local field potentials, and electroencephalograms)~\cite{xue2016minimum,klaus2011statistical}. 
The following dynamics describe these systems 
\[
\Delta^\alpha x(t+1) = Ax(t)+Bu(t), \ t=0,1,\ldots,
\]
where $\alpha = [\alpha_1,\ldots,\alpha_n]^\top \in\mathbb{R}_{+}^n$ are the fractional exponents, and where $\Delta^\alpha x_k = \sum_{j=0}^k D(\alpha,j)x_{k-j}$ is the fractional derivative given by 
\[
\begin{array}{c}
D(\alpha,j) = \text{diag}(\psi(\alpha_1,j),\ldots,\psi(\alpha_n,j)), \text{and}, \\
    \psi(\alpha,j) = \displaystyle\frac{\Gamma(j-\alpha)}{\Gamma(-\alpha)\Gamma(j+1)} = \prod_{\ell=1}^j\left(\frac{j-\ell-\alpha}{j+1-\ell}\right)
\end{array}
\]
where $\Gamma(\cdot)$ denotes the gamma function defined as 
$
    \Gamma(z) = \int_0^\infty t^{z-1}e^{-t}dt
$ 
for $z\in\mathbb{C}$. It is worth noticing that if $\alpha=\mathbf{1}_n$ (i.e., the $n$-vector of ones), then we obtain the description of an LTI. Thus, fractional-order systems are inherently nonlinear with infinite memory but require only compact parametric descriptions. 
Structural  observability properties of these systems, as well as sensor placement, has been addressed in
~\cite{pequito2015minimumBrain}.

\subsection{Switching systems}

In what follows, we consider two classes of switching systems: (\emph{i}) temporal networks, i.e., the same structural pattern across switches but with different realizations over time; and  (\emph{ii}) linear time-invariant switching systems (or their upper classes such as descriptor LTI).

\subsubsection{Temporal networks}\label{sub:temporal}

In recent years, a particular focus from the networks science community has been on temporal networks~\cite{holme2012temporal}. In this context, structural controllability of temporal networks considers those modeled by linear time-varying systems whose structural pattern remains unchanged, and their realization may vary over time~\cite{posfai2014structural}. In~\cite{hou2016structural}, the authors consider the case where there is a finite number of possible switches (i.e., a finite number of realizations of a given structural pattern), yet it can be repeated multiple times. In~\cite{yao2017structural}, the authors propose using a switching controller to increase the dimension of the structural controllable subspace.

\subsubsection{Linear time-invariant switching systems}

Conceptually, we can see a linear time-invariant switching system (LTIS) as a set of LTI systems, where each element of the set is called a \emph{mode}, together with a set of discrete events that cause the system to switch between modes.  
Subsequently, an LTIS may be described as follows: 
\begin{equation}\label{eq:slcs}
	\dot x(t) = A_{\sigma(t)}x(t)+B_{\sigma(t)}u(t),
\end{equation}
where $\sigma:\mathbb R^+\to \mathbb M=\{1,\hdots,m\}$ is a piecewise switching signal, that only switches once in a given dwell-time, $x(t)\in\mathbb R^n$ the state of the system, and $u(t)\in \mathbb R^p$ is a piecewise continuous input signal. 
As we may expect, finding a set of sparsest input matrices $\{B_{\sigma(t)}\}$ that ensures each mode of the system to be controllable is an NP-complete problem~\cite{ramos2018robust}. 

The necessary and sufficient conditions to ensure structural controllability of linear time-invariant switching systems have been proposed in~\cite{ramos2013model,liu2013structural,liua2013graph}, with possible robustness characterization~\cite{ramos2015analysis,ramos2013structural}. 
Afterward, the actuator placement problem was addressed in~\cite{pequito2017structural}. 

In~\cite{boukhobza2011observability,boukhobza2012sensor} the authors proposed conditions for sensor placement in linear time-invariant switching systems with unknown inputs. Later, they generalized these conditions to linear time-invariant switching descriptor systems~\cite{boukhobza2013discrete,gracy2018input}.

\subsection{Linear time-varying systems}\label{sub:time_varying}

A linear time-varying systems (LTV) can be seen as the system in~\eqref{eq:lti}, where the matrices $A$, $B$, $C$ and $D$ are also time-dependent, i.e., vary with time. 

In~\cite{lichiardopol2007linear}, the authors proposed conditions to evaluate the structural controllability of linear time-varying systems. In~\cite{hartung2013necessary}, the authors extend these results and compare their implications in different controllability contexts. 
In~\cite{gracy2017structural}, the authors present topological conditions required to assure FDI for linear time-varying systems with unknown inputs, which allow them to retrieve both the initial state and the unknown inputs over long time windows.

\subsection{Petri nets}\label{sub:petri}

Petri nets consist of a framework that allows the modeling and analysis of discrete event dynamic systems. 
Petri nets can be represented as bipartite graphs, where the state variables are called places, and the transformations on the states are referred to as transitions.  
Such places transitions are connected through pre-incidence (i.e., inputs) and post-incidence (i.e., outputs), under possible constraints indicating the resources required. 
Due to the graphical nature of these networks, it is possible to leverage the notion of structural observability to retrieve the places (i.e., state variables) under known transitioning models~\cite{silva2013half}. 
To assess structural observability, we must consider a suitable transformation of the original graph~\cite{silva2011fluidization,mahulea2010observability}.

\subsection{Hybrid systems}\label{sub:hybrid}

A hybrid system is a dynamical system that manifests both continuous and discrete dynamic behavior. A special case of hybrid systems is when the continuous dynamics is given by an LTI system, which yields the class of linear time-invariant switching systems. 

In~\cite{ramos2013structural}, the authors introduce a tool for the design and verification of structural controllability for hybrid systems. Later they considered it in the context of the 
analysis and design of electric power grids with robustness guarantees on the link failures~\cite{ramos2015analysis}. 


\subsection{Nonlinear systems}\label{sub:nonlinear}

The authors, in~\cite{stefani1985local}, revealed a local controllability condition for nonlinear systems. 
Nevertheless, there is limited research about structural controllability for nonlinear systems. 
The first work in the line of structural controllability of such systems is the one in~\cite{fradellos1977structural}, where changes concerning controllability or uncontrollability behavior of perturbed linear and non-linear systems are examined.  

In~\cite{qiang2010some}, the conception of structural controllability is extended to nonlinear system utilizing Lie algebra theory. This extension is then used to analyze the structural properties of nonlinear system. 
In~\cite{ma2010structural}, the author proposes to assess the structural controllability of the nonlinear system through the system transfer function. 
In~\cite{zanudo2017structure}, the authors propose to use feedback vertex sets in relation to structural properties to assess the controllability of nonlinear systems. 
In~\cite{angulo2019structural} the authors leverage structural properties to enable nonlinear assessment of controllability and observability properties. 
 In~\cite{kawano2019structural}, the authors construct necessary conditions for the structural controllability and observability of complex nonlinear networks. 
These conditions, which are based on refined notions of structural controllability and observability, can be used for networks governed by nonlinear balance equations to develop a systematic actuator/sensor placement.

In~\cite{van2018zero}, the author explores the controllability conditions of nonlinear systems performing its linearization, and in the same lines
\cite{staroswiecki2007structural,staroswiecki2007observability} the authors suggest to explore FDI settings.

\section{Variations on Structural Systems Theory}

As previously emphasized, structural systems theory deals with the structural pattern of matrices. It is assumed that the realization of the parameters takes place on an infinite field (e.g., the reals) and that these parameters are independent of each other. In what follows, we provide a brief description of results built upon structural systems theory to obtain methodologies to hand cases where such assumptions might be waived.

\subsection{Positive Systems}

A discrete-time linear time-invariant system is a \emph{positive system} if for any initial condition and any nonnegative input sequence, all the entries of the state vector remain positive over time.  

In~\cite{commault2004simple}, the authors addressed the reachability of discrete-time linear time-invariant systems to assess when the state can lie on the positive octant that has broad applications in practice~\cite{rantzer2018tutorial}. This problem also motivated the study presented in~\cite{ruf2018herdable,she2020characterizing}, where the authors address the input selection to achieve the reachability property. In~\cite{lindmark2016positive}, the authors build up on structural systems to explore the parametric dependencies that lead to the controllability of positive systems. Along the same lines, the authors in~\cite{bru2005monomial} explore related digraph properties to study the controllability of positive systems.

\subsection{Sign Systems}

As an extension to the use of structural systems theory to positive systems, the authors in~\cite{hartung2014sign} leverage the former systems to render necessary conditions for signed systems' properties, specifically to attain sign stabilizability.

\subsection{Parameter-dependent Structural Systems}\label{sec:parameterDepSS}

In~\cite{whalen2015observability,whalen2016effects,menara2018structural,mousavi2017structural}, the authors examine the structural controllability criterion when the dependency of the parameters is given by the symmetry of the system's autonomous matrix. For instance, in~\cite{menara2018structural}, the authors show that (symmetric) structural controllability can be assessed by graph theoretic elements similar to those previously proposed to verify structural controllability. In~\cite{romero2018actuator}, the authors used the latter criterion to addressed the actuator placement problem in this context under possible cost constraints.

In~\cite{zhang2019structural}, the authors investigate conditions when subsystems satisfy fractional parametrizations. 
In~\cite{liu2019graphical}, the authors propose a graphical characterization to attain structural controllability when arbitrary linear dependencies exist between the system's dynamics parameters. 

Lastly, it is worth mentioning that in \cite{murota2009matrices,murota2012systems} the parametric dependencies have been accounted for, using the notion of \emph{mixed matrices}, where the entries could be zero/nonzero or fixed constant, often capturing the network dependencies (or, generally speaking, losses of degrees of freedom). To address the actuator placement problem (among other problems), the author deploys matroid theory and some of the schemes to approximate the solutions.

\subsection{Structural theory on finite fields}\label{sec:finiteField}

In~\cite{feng2008structural,sundaram2012structural}, the authors propose to assess structural controllability properties when the parameters are taken to be independent but take values on a finite field. In this context, properties are no longer valid generically but rather with a certain likelihood.  
In~\cite{yuan2016structural}, the authors provide conditions for the analysis and design of such systems in the frequency domain.

\subsection{Strong Structural Theory}

In contrast with structural theory, the (classical) strong structural theory~\cite{mayeda1979strong} seeks to guarantee properties for any set of parameters considered for the realization of the nonzero entries of the structural pattern except when the parameters are zero. 
Notwithstanding, it is possible to extend to scenarios where some of the entries could be either zero, nonzero (i.e., a real scalar different from zero) and possibly a real scalar that could be either a zero or a nonzero~\cite{popli2019selective,jia2020unifying}. 
It is worth noticing that several of the problems discussed in this review can be posed in the strong structural controllability scenario, but their computational complexity often changes. 
For instance, the strong structural controllability of the sparsest minimal controllability problem is NP-hard~\cite{trefois2015zero}. This result contrast with the existing polynomial solution when the goal is to attain structural controllability. 

\subsection{Bond-graphs}

The concept of bond graphs~\cite{paynter1961analysis} seeks to describe the dynamic behavior of physical systems based on energy and energy exchange~\cite{thoma2016introduction}. The \emph{basic units} can be seen as concepts and/or objects, enabling object-oriented physical systems' modeling. Bond graphs are labeled directed graphs, in which the vertices represent basic units, and the edges represent an ideal energy connection between them, and they are referred to as bonds. As such, each basic unit could be seen, in particular, as a linear time-invariant system that interconnected through the others using bonds. Structural systems theory has been leveraged to assess several system properties of these systems and to address similar problems as those overviewed in this survey, seeking to design the systems to attain such properties -- see, for example,~\cite{sueur1989structural,sueur1991bond,alem2014bond}.

\section{Applications}


This section describes a diversity of applications that explores the concepts of structural systems.

\begin{figure}[!ht]
    \centering
    \includegraphics[width=0.85\textwidth]{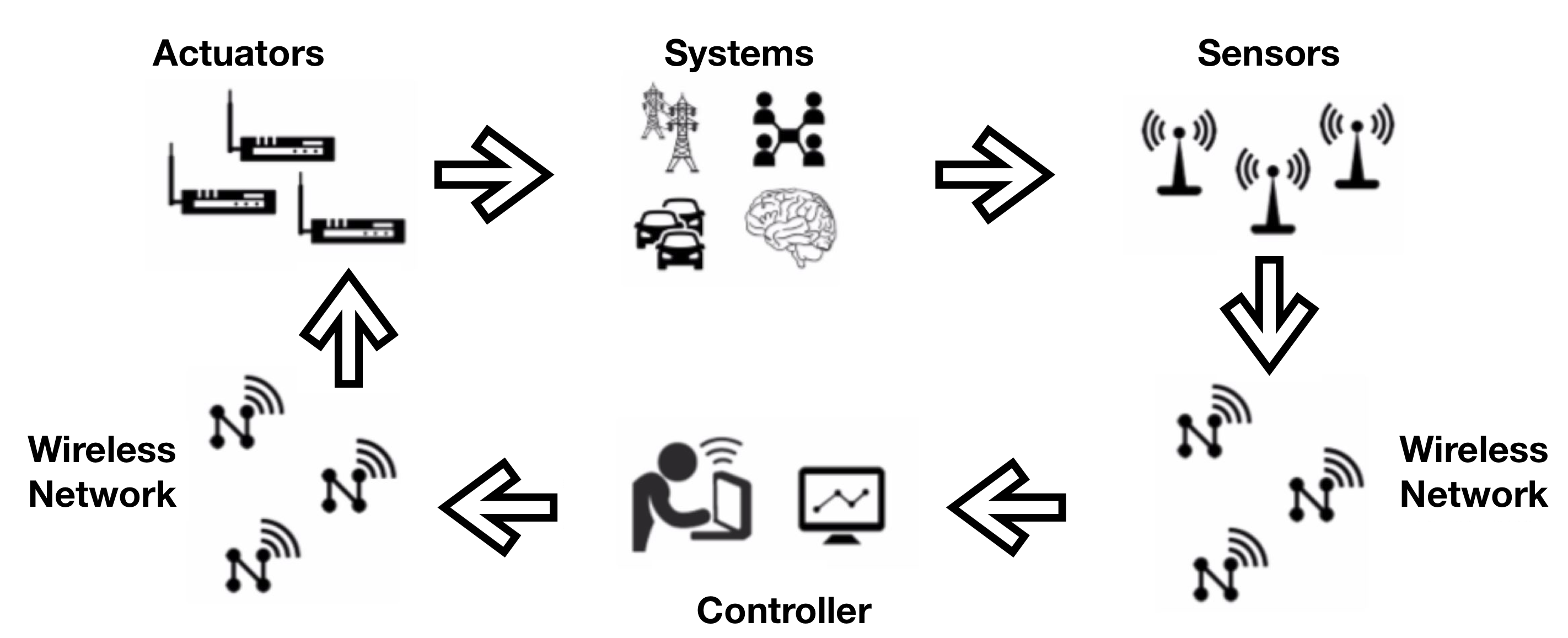}
    \caption{Applications of structural systems.}
    \label{fig:wireless_net}
\end{figure}

\subsection{Distributed/Decentralized Estimation and Optimization}

In~\cite{sundaram2008distributed,doostmohammadian2013genericity,doostmohammadian2014graph}, the authors use the notion of structural observability to ensure that each node in the network is capable of retrieving the neighbors' states towards computing the solution to an optimization problem in a distributed fashion.

In~\cite{alexandru2016decentralized,alexandru2017limited}, the authors propose the design of the sensor network topology to guarantee that every sensor can retrieve the state of the entire system under possibly link failures.

In~\cite{khan2011coordinated}, the authors investigate generic observability properties to infer the stability of distributed estimation schemes, possibly under distributed settings~\cite{doostmohammadian2018structural}.


\subsection{Wireless networks}\label{sub:wir_net}

In~\cite{sundaram2010wireless,pajic2011wireless,pajic2013topological},  the authors consider sensor networks as an extension to the state space,  connected through feedback with the plant. 
Due to the freedom in the design of the sensor network, it is possible to guarantee enough redundancy, in terms of paths between the states of the plant and the sensors, such that proper monitoring is secured (see the section Security and Resilience, Section~\ref{subsub:security}, for more details). 
In~\cite{pajic2011wireless} the authors aimed to assess the controllability and stability properties of these joint plant-sensor network systems. 
Some privacy properties have also been considered~\cite{pequito2015smart}. 
Additionally, extensions of this setting considered the scenario of potential communication link failures~\cite{sundaram2012control}. In~\cite{martinez2010communication}, the authors explore how the transmission sequence, through a networked system over wireless networks, should be designed to ensure controllability and observability properties. In~\cite{kruzick2017structurally}, the authors address the observability problem when backbone nodes (e.g., routers) are considered in the context of sensor networks.

Alternatively, if the sensors include a queuing capability (i.e., \textsf{WirelessHART}, see details in~\cite{d2018modeling}), then the authors in~\cite{d2016resilient} were able to derive necessary conditions for observability and stabilizability of the network, in the context of state feedback.

\subsection{Networked Control Systems}

A Networked Control System is a control system where the control loops are closed through a communication network.

In~\cite{sauter2006decentralized}, the authors proposed an FDI scheme (see Section FDI) for networked control systems. 
In~\cite{doostmohammadian2019cyber}, the authors leverage the structural system's results to propose a cyber-social system framework. 
Lastly, in~\cite{zhang2019structural}, the authors explore conditions when subsystems satisfy fractional parametrizations. 

\subsection{Network Neuroscience}

The works in~\cite{gu2015controllability,tu2018warnings,pasqualetti2019re,tang2018colloquium} assess the structural controllability aspects of brain networks. Also, a work-related to brain networks, in the context of using discrete-time fractional-order systems, is the one in~\cite{pequito2015minimum} (see the section of Nonlinear models, Section~\ref{sub:nonlinear})

\subsection{Multi-agents}

In the context of multi-agents, due to the finite memory and number of tasks required to perform, it turns out that their dynamics are modeled using finite fields. 
In~\cite{sundaram2012structural,sundaram2013control}, the authors leveraged the notions of structural controllability and extended the dynamics in the context of finite fields realizations -- see also Section~\ref{sec:finiteField}. 

Alternatively, under (regular) structural systems theory (i.e., infinite field realizations), a multitude of approaches have been proposed. 
For instance, in~\cite{zamani2009structural}, the authors address the structural controllability of multi-agent networks driven by a single agent. In~\cite{ouyang2018controllability}, the authors explored the minimum number of leaders that ensure structural controllability (i.e., the \emph{leader selection problem}), which may be recast as an input selection problem. 
The authors in~\cite{partovi2010structural} study the structural controllability of high-order dynamic multi-agent systems, and in~\cite{guan2017structural,liu2013graph} results for switching topologies are presented.  
In~\cite{liu2017scheduling}, the authors addressed the actuator and sensor placement under possible cost constraints for multi-agent bilinear systems (see Section~\ref{sub:bilinear} for more bilinear systems related work). 
In~\cite{pequito2014optimal}, the authors address the design of a communication topology dynamics that should be observable from each agent, and that achieves minimal overall transmission cost. 

In the context of (fully) distributed leader selection, 
in~\cite{pequito2015distributed,tsiamis2017distributed}, the authors propose a two-level approach. 
First, the agents determine with which agents they need to interact with to ensure structural controllability. 
Then they pick weights locally that ensure controllability of the overall network.
In the same spirit, the authors in~\cite{mehrabadi2019structural} leverage structural controllability properties to propose a parametrization technique to attain consensus from a centralized perspective, with a collection of multiple leaders.

\subsubsection{Consensus/agreement protocols}

In~\cite{goldin2013weight}, the authors characterize generic controllability properties for dynamics that implement consensus algorithms. It is worth emphasizing that the direct application of structural systems results does not ensure the desired assessment of the system properties as in these dynamics the diagonal entries depends on the remaining row entries, thus violating the assumption that all nonzero parameters are independent. As an alternative, one can considered the settings discussed in~\ref{sec:parameterDepSS}.

\subsection{Power grid}\label{sub:power_grid}

The work in~\cite{bhela2017enhancing} ensures that rank conditions are achieved in the context of algebraic differential equations, as part of the power flow optimization.  In~\cite{bhela2017enhancing,bhela2017power}, the authors  assess the  identifiability of the linearized power grid system by leveraging structural systems theory (see also Section~\ref{subsub:Struct_ident}), as well as the vulnerability to attain this~\cite{luo2018structural}. 
In~\cite{xiaoyu2012investigation}, the authors assess the structural controllability of electrical networks using rational function matrices.

\subsection{Medical}\label{sub:medical}

Structural identifiability was proposed to be used in the context of a model for dialysis~\cite{canto2009structural}. 
In~\cite{pequito2015minimum}, the authors proposed to assess the minimum number of electrodes required for monitoring electroencephalographic data (i.e., sensor placement in the brain).

\subsection{Network Coding}

In~\cite{campobello2009distributed}, the authors leveraged structural systems to impose conditions for network coding in the context of state-space representation in the spirit of~\cite{koetter2003algebraic}. In~\cite{sundaram2009linear}, some of these ideas are used in the context of sensor networks.

\subsection{Science Exploration Tools}\label{sub:tools}


In~\cite{liu2016control}, the authors overview some of the different applications of structural control as a tool to unveil features in the context of network science. 
This application has attracted interest since 2011 when the paper~\cite{liu2011controllability} feature at the cover of Nature. In the latter, the authors explored the actuation placement problem and provided evidence of its correlation with the degree distribution when nodal dynamics were not considered~\cite{cowan2012nodal}.
More recently,  in~\cite{liu2012control}, the notion of control centrality is introduced.

Among the possible applications, structural systems theory enables the characterization of networks in classes. 
For instance, in \cite{ruths2014control}, the characterization is upon the partition of dilations in a network. This characterization can also be used to construct models to attain such characterizations~\cite{campbell2015topological}. 
More recently, in~\cite{pequito2017trade}, the authors introduce the notion of \emph{actuation spectrum}, which captures the trade-offs between the minimum number of state variables required to attain structural controllability in a given number of time-steps. 
Lastly, in~\cite{ramos2020generating}, the authors provide evidence that several of the generative models and centrality measures fail to capture the actuation spectrum of real networks. Therefore, they proposed a novel generative model that builds upon the notion of structural time-to-control communities.

\subsection{Other applications}\label{sub:other_app}

FDI applications are considered in~\cite{verde2007monitorability} and~\cite{veldman2015towards} application for monitoring the gas turbine and water distribution networks, respectively.

\subsection{Software Routines}\label{sec:software}

Besides the code available by the different authors at their personal websites and file exchange platforms (e.g., Mathworks), there are MATLAB toolboxes such as~\cite{geisel2019matlab},  C$++$~\cite{martinez2007lisa}, and Modelica and Python~\cite{perera2015structural}.


\section{Conclusions and Future Research Directions}\label{sec:fut_conc}


In this document, we rendered a survey of the advancements achieved since the latest survey conducted by Dion \emph{et al.}~\cite{dion2003generic}, in the scope of structural systems theory. 
We considered all the papers written in English that use or leverage structural systems theory, including their direct variants and applications.
Specifically, we analyzed each peer-reviewed paper in the press that cited either Lin's paper~\cite{lin1974structural}, or Dion \emph{et al.}~\cite{dion2003generic}, if it was available from 2003 to the end of March 2020. 
The conference publications were only considered whenever they are not part of a published journal version. 
Analogously, we did not overview dissertations because they comprise papers that were cited. 
Also, the presented study is not an in-depth overview of all the results in the area. 
Instead, we provided a (possibly bias) selection of results that may help future advancements in the context of structural systems theory and its use to show \emph{non-}structural system properties. 
For similar reasons, we tried to report the findings in the reported period without validity and usefulness judgments. 

The first section of this paper introduces the area to a broad audience. Next, we introduce essential notions to keep the document self-contained, and we motivate the use of structural systems. 
Subsequently, we presented the scientific advancements in the area of LTI systems, followed by the findings in other classes and variations of structural systems.

As hinted from the above overview, it is clear that some topics were explored in greater depth than others. 
Specifically, there is still a lack of necessary and sufficient conditions to attain structural systems properties, enabling us to design the systems to attain such properties. 
Besides, it would be interesting to extend the results obtained for linear time-invariant systems to the other subclasses covered in Section~\ref{sec:other_sublasses}. Furthermore, it is essential to explore trade-offs between different system theoretic properties and understand how to design systems such that such trade-offs are attained.

Secondly, it is only natural to understand the computational complexity of assessing different structural conditions, and the complexity of designing systems to attain such properties. 
Notice that even though some design problems are NP-hard, one should determine the most significant subclasses of problems that permit to perform the design using polynomial algorithms. 
At the core of this research, the quest is the symbiotic relationship between the algorithms used and the system's structure, which allows us to craft current algorithms to perform efficiently and optimally in some domain of interest. 
In practice, there is also a demand to improve the computational efficiency of the algorithms already used to assess structural properties by either considering approximations as a good starting point to the optimal algorithms or developing distributed versions of those algorithms that can be implemented in parallel platforms.

 Furthermore, it is worth mentioning that in \cite{murota2009matrices,murota2012systems} the parametric dependencies have been accounted for, using the notion of \emph{mixed matrices}, where the entries could be zero/nonzero or fixed constant, often capturing the network dependencies (or, generally speaking, losses of degrees of freedom). In these cases,  matroid theory and its algorithms can be used to approximate the solutions to design problems (e.g., actuator placement). In particular, even when optimality cannot be guaranteed in several of the design problems, it is often the case that some suboptimality guarantees may be ensured and improved by considering the system's structure. 
For instance, some structural systems properties are \emph{submodular}, for which efficient greedy algorithms are available~\cite{clark2015submodularity,bach2011learning}.
For example, structural controllability problems can be posed as matroid optimization problems that can be solved exactly under certain assumptions~\cite{rocha2014output,clark2017input}. 
Such developments should be complemented with recent research that unveiled new insights on additional properties (e.g., submodular ratio and curvature);  thus, tightening the suboptimality guarantees~\cite{iyer2013curvature,bian2017guarantees,gupta2018approximate}.

Another direction is leveraging structural systems theory as a tool to aid in finding a solution to optimization algorithms. On the one hand, we can consider algorithms that first aim to attain structural systems properties (e.g., structural controllability/observability/stability), which ascertain feasibility of the solutions for almost all sets of parameters. 
Then, we determine the set of parameters that minimize/maximize a desirable objective. 
For instance, consider the following three possible applications:  
(i) in~\cite{becker2020network}, the authors use structural systems to ensure that for almost all set of parameters stability and controllability would be guaranteed to determine a set of perturbations in the dynamics that improve the controllability energy (see Section~\ref{sec:systemDesignsC}); 
(ii) in~\cite{pequito2018analysis}, the authors guarantee generic stabilizability properties for the decentralized control, and then an iterative procedure is considered to find a set of parameters that stabilize the plant -- see Remark~\ref{structNonStruc}; and 
(iii) in~\cite{pequito2015distributed},  towards determining a fully distributed leader selection to attain controllability. They decompose the problem into two. The first determines the leaders that attain structural controllability. The second determines a set of parameters for the local interactions that ensure (non-structural) controllability of the network.

On the other hand, we could envision a \emph{structural-convex optimization} framework, where several discrete mathematics algorithms could be intertwined with the convex optimization tools already in use. In particular, when considering convex relaxations that often involve performing a set of operations that lead to a description of the optimization problem where the only nonlinear constraint is that of a rank constraint on a matrix of interest~\cite{boyd2004convex}. In this context, we suggest the use the generic rank instead of the rank constraint that is often dropped from the optimization problem such that we can guarantee an upper bound on the rank of matrices with some structure (e.g., see Remark~\ref{Toeplitz}). Additionally, it might be possible to leverage the structural similarity between the rank of a transfer function and the Schur complement to perform some of the algebraic transformations and ensure some graph-theoretical properties on the generic rank. Specifically, consider the matrix $M(s)=\left[\begin{array}{cc} A-sI & B\\ C &0\end{array}\right]$ for which the following holds: $\text{rank } C(sI-A)^{-1}B=\text{rank } M(s)-n$. Now, notice that the $\text{grank } C(sI-A)^{-1}B$ equals to the number of vertex-disjoint paths form the inputs to the outputs in the system digraph $\mathcal G(\bar A,\bar B,\bar C)$ (see~\cite{van1991graph}). Thus, we can potentially design the structure of the matrices $(\bar A,\bar B,\bar C)$ such that almost all realizations would attain a desirable rank, after which usual convex optimization tools could be used.









\bibliography{mybibfile}

\end{document}